\documentclass[final,12pt]{elsarticle}
\journal{Stochastic Processes and their Applications}

\usepackage{amsthm,amsmath,amsfonts,amssymb}
\usepackage[numbers]{natbib}
\usepackage{lineno,hyperref}
\usepackage{appendix}

\numberwithin{equation}{section}

\theoremstyle{plain}
\def\R{{\mathbb R}}
\def\N{{\mathbb N}}

\def\rank{{\text{\bf{rank}}}}
\theoremstyle{plain}

\newtheorem{lem}{Lemma}
\newtheorem{thm}{Theorem}
\newtheorem{asp}{Assumption}
\newtheorem{cor}{Corollary}

\newtheorem{pro}{Property}
\newtheorem{exam}{Example}

\newtheorem{remark}{Remark}
\newtheorem*{notation}{Notation}

\def\R{{\mathbb R}}
\def\N{{\mathbb N}}

\def\rank{{\text{\bf{rank}}}}

\begin{document}

\begin{frontmatter}

\title{A NOVEL COMPUTATIONAL FRAMEWORK FOR THE EXPECTED NUMBER OF REAL ROOTS OF STOCHASTIC FUNCTIONS ON A GIVEN INTERVAL}

\author{Yi Xu}
\address{Mathematics Department, Southeast University, 2 Sipailou, Nanjing, Jiangsu Province, 210096 P. R. China}
\ead{yi.xu1983@hotmail.com}

\begin{abstract}
We propose a new computational framework for the expected number of real roots of a stochastic function on a given interval. The classical Kac--Rice formula requires the joint density of the function and its derivative, which is often intractable. Our approach avoids this requirement entirely by introducing a cumulative expectation function $\varphi_{G,\rho}(t)$. Through analysis of its absolute continuity and differential structure, we derive two complementary computational schemes: one expresses the expectation as a derivative of a variable-domain integral under weak conditions; the other yields an explicit integral representation without joint densities or variable‑domain differentiation.

We illustrate the method in detail for linear stochastic functions, obtaining explicit formulas for Gaussian and uniform distributions, together with several new analytical results. The framework substantially broadens the scope of problems amenable to rigorous analysis and provides a powerful tool for applications in stochastic analysis and beyond.
\end{abstract}


\end{frontmatter}
\section{Introduction}
Consider the equation with variable $t \in \R$ and parameter $x \in \R^n$:
\[
G(t, x) = 0, \quad t \in I = [l, u],
\]
where $G(t, x) : I \times \R^n \rightarrow \R$, and $l, u \in \R$.
For a given $x \in \R^n$, define the number of real roots of the above equation within $I$ as
\[
N(x, G, I) = \#\{ t \in I \mid G(t, x) = 0 \}.
\]

If $x$ is a random variable with density function $\rho(x)$, we aim to compute the expectation of $N(x, G, I)$, i.e.,
\begin{equation}\label{eq-0}
  \mathbb{E}N(G, \rho, I) = \int_{\R^n} N(x, G, I) \, \rho(x)  dx.
\end{equation}

This problem represents a classical topic at the intersection of stochastic analysis and equation theory, with conclusions widely applied in engineering optimization, statistical inference, physical modeling, and numerous other fields. Since Littlewood's pioneering investigation in 1930s into the expected number of real roots under discrete distributions\cite{littlewood1938, littlewood1939, littlewood1943}, related research has evolved over nearly a century. Early results were confined to low-dimensional variables and special distribution scenarios. It was not until the introduction of the Kac--Rice formula that a unified framework for expressing this expectation was established \cite{kac1943, kac1948, rice1945}. This formula provides an integral expression for \eqref{eq-0} through the joint density function $\rho_{G,G'}(x)$ of $G(t,x)$ and $G'(t,x)$. However, it also exhibits significant limitations, which introduce two fundamental difficulties:

(1) Computational Difficulty: Except in a few special cases (e.g., when $G$ and $G'$ are independent), obtaining the joint density function itself is an extremely complex task. It essentially involves performing a multiple integral over a variable domain with two parameters followed by a secondary differentiation process, making the derivation of its explicit expression exceptionally challenging in general.

(2) Theoretical Overhead: Even disregarding computational concerns, the application of the Kac--Rice formula requires this joint density function to satisfy a series of additional conditions. Verifying these conditions without an explicit expression for the joint density is equally difficult. 

As noted by Adler and Taylor \cite{adler2007}, ``it turns out to be remarkably difficult to compute the integral in Rice's formula unless $G$ is either Gaussian or a function of Gaussian processes.'' Subsequent extensions to higher‑dimensional parameter sets \cite{adler2007, azaïs2009, berzin2022} have enriched the theory, yet the reliance on the joint density remains an intrinsic obstacle.

To overcome this limitation, we propose a framework that avoids joint densities entirely and relaxes the required assumptions. The key idea is to study the cumulative expectation function
\[
\varphi_{G,\rho}(t) = \mathbb{E}N(G,\rho,[l,t]), \qquad t\in [l,u],
\]
so that $\varphi_{G,\rho}(u)$ is the quantity of interest. Our approach proceeds in two steps:
\begin{itemize}
\item[(1)] Establish absolute continuity of $\varphi_{G,\rho}(t)$, ensuring that $$\varphi_{G,\rho}(u)=\varphi_{G,\rho}(l)+\int_l^u \varphi_{G,\rho}'(t)dt.$$
\item[(2)] Use the Coarea Formula to transform $\varphi_{G,\rho}'(t)$ into concrete integral representations.
\end{itemize}

The main contributions are two complementary computational schemes:

\begin{itemize}
\item[] \textbf{Scheme 1} (Theorem \ref{thm-2}, Corollary \ref{cor-main1}): Under weak assumptions,
\[
\varphi'_{G,\rho}(t^*)=\left(\int_{\Theta(t^*,\varepsilon)} \rho(x)dx\right )'\Big|_{\varepsilon=0},
\]
where $\Theta(t,\varepsilon) =\{ x \in \R^n \mid G(t,x)G(t+\varepsilon,x) \leq 0 \}.$ This reduces the problem to differentiating a variable‑domain integral, a much simpler task than obtaining joint densities.
\item[] \textbf{Scheme 2} (Theorem \ref{thm-3}, Corollary \ref{cor-main2}): Under slightly stronger but verifiable conditions,
\[
\varphi'_{G,\rho}(t^*)= \int_{L(G(t, x), \{t^*\})} \rho(x) \frac{|G'(t^*, x)|}{\|\nabla G(t^*, x)\|}  d\mathbb{H}^{n-1},
\]
yielding an explicit integral formula free of joint densities and variable‑domain differentiation.
\end{itemize}

As an illustration, we treat in detail the linear case $G(t,x)=f(t)^\top x - f_0(t)$. Explicit formulas are derived for Gaussian and uniform distributions, yielding several new analytical results that are difficult to obtain via traditional methods.

The paper is organized as follows. Section 2 sets notation and basic properties. Section 3 collects lemmas from geometric measure theory (proofs are in the Appendix). Section 4 establishes absolute continuity of $\varphi_{G,\rho}(t)$. Section 5 presents the two main theorems. Section 6 analyzes the linear case with examples. Section 7 sketches an extension to variable intervals. Section 8 concludes.

\section{Basic Notation and Properties}
Throughout this paper, we maintain the following core assumption:
\begin{asp}\label{ass-0}
\begin{itemize}
    \item[(1)] The density function $\rho(x) : \R^n \rightarrow \R$ is integrable with respect to the Hausdorff measure on $\R^n$;
    \item[(2)] The expectation satisfies $\mathbb{E}N(G, \rho, I) < +\infty$.
\end{itemize}
\end{asp}

We begin by introducing the notation that will be used extensively in the subsequent analysis:
\begin{notation}
\begin{itemize}
    \item[(1)] $\N = \{1, 2, \cdots\}$ denotes the set of natural numbers; $\R^n$ is the $n$-dimensional real vector space; $\R^{m \times n}$ is the space of $m \times n$ real matrices; $0_{m \times n}$ denotes the $m \times n$ zero matrix; $0_n$ denotes the $n$-dimensional zero column vector.
        \item[(2)] $\|x\|=\sqrt{x^Tx}$ denotes the Euclidean norm; $|\Omega|$ denotes the Lebesgue measure of $\Omega$; $\mathrm{cl}(\Omega)$ denotes the closure of $\Omega$; $B(x,\delta)$ denotes the hyperball centered at $x$ with radius $\delta$.
    \item[(3)] For $J \subseteq A \subseteq \R^k$, and functions $P(s, x) : A \times \R^n \rightarrow \R^m$, $Q(s, x) : A \times \R^n \rightarrow \R^k$, we define:
          \[
          L(P(s, x), J) = \bigcup_{s \in J} \{ x \in \R^n \mid P(s, x) = 0_m \},
          \]
          \[\eta(P(s,x), J, x)=\{s\in J \mid P(s,x)=0_m\},\]
          \[
          \Xi(P(s,x), Q(s,x), J)=\{x\in \R^n \mid \emptyset\neq \eta(P(s,x),J,x)\subseteq \eta(Q(s,x),J,x)\}.
          \]
    \item[(4)] $\Lambda_i(G, I) = \{ x \in \R^n \mid N(x, G, I) \geq i \}$ for $i \in \mathbb{N}$.
    \item[(5)] $(P'(s, x))_{i,j} = \displaystyle\frac{\partial P_j(s, x)}{\partial s_i}$ and $(\nabla P(s, x))_{i,j} = \dfrac{\partial P_j(s, x)}{\partial x_i}$.
    \item[(6)] The indicator function is defined as $\mathbb{I}(x, \Omega) = \begin{cases} 
        1, & x \in \Omega \\
        0, & \text{otherwise}
    \end{cases}$.
    \item[(7)] If $\Omega$ is an $m$-dimensional Hausdorff measurable set, then $\displaystyle\int_{\Omega} \rho(x)  d\mathbb{H}^m$ denotes the integral of $\rho(x)$ over $\Omega$ with respect to the $m$-dimensional Hausdorff measure. In particular, if $\Omega$ is an $n$-dimensional set, then
          \[
          \int_{\Omega} \rho(x)  d\mathbb{H}^n = \int_{\Omega} \rho(x)  dx.
          \]

\end{itemize}
\end{notation}

Based on the above notation, we now establish several fundamental properties.

\begin{pro}\label{pro-1}
For any subinterval $I_i \subseteq I$, the following properties hold:
\begin{itemize}
    \item[(1)] $0 \leq \mathbb{E}N(G,\rho,I_i) < +\infty$, $\mathbb{E}N(G,\rho,\emptyset)=0$, $L(G(t,x),I) = \Lambda_1(G,I)$, and
          \[
          \mathbb{E}N(G,\rho,I) = \int_{L(G(t,x),I)} N(x, G, I) \rho(x)  dx = \int_{\Lambda_1(G,I)} N(x, G, I) \rho(x)  dx.
          \]
    \item[(2)] $\mathbb{E}N(G,\rho,I_1 \cup I_2) = \mathbb{E}N(G,\rho,I_1) + \mathbb{E}N(G,\rho,I_2) - \mathbb{E}N(G,\rho,I_1 \cap I_2)$.
    \item[(3)] If $I_i \cap I_j = \emptyset$,  $\forall i \ne j$, then
          \[
          \mathbb{E}N\left(G,\rho,\bigcup_{i=1}^{\infty} I_i\right) = \sum_{i=1}^{\infty} \mathbb{E}N(G,\rho,I_i).
          \]
    \item[(4)] If $I_1 \subseteq I_2$, then $\mathbb{E}N(G,\rho,I_1) \leq \mathbb{E}N(G,\rho,I_2)$.
    \item[(5)] $\Lambda_{i+1}(G,I) \subseteq \Lambda_{i}(G,I),$  $\forall i \in \mathbb{N}$.
    \item[(6)] If $I_1 \subseteq I_2$, then $\Lambda_i(G,I_1) \subseteq \Lambda_i(G,I_2)$.
    \item[(7)] $\mathbb{E}N(G,\rho,I) = \displaystyle\sum_{i=1}^{\infty} \int_{\Lambda_i(G,I)} \rho(x)  dx$.
    \item[(8)] For $P(s, x) : A \times \R^n \rightarrow \R^m$, $Q(s, x) : A \times \R^n \rightarrow \R^k$, and  $\forall J\subseteq A\subseteq  \R^m$,
          \[
          \Xi(P(s,x),Q(s,x),J)\subseteq L(\{P(s,x),Q(s,x)\},J).
          \]
\end{itemize}
\end{pro}

\proof
(1) All properties except $\mathbb{E}N(G,\rho,I_i) < +\infty$ are evident from the definitions.
Since $0 \leq N(x, G, I_i) \leq N(x, G, I)$ and
\[
\mathbb{E}N(G,\rho,I) = \int_{\R^n} N(x, G, I) \rho(x)  dx < +\infty,
\]
it follows that $\mathbb{E}N(G,\rho,I_i)$ exists and is finite.

(2) Follows directly from the identity:
\[
N(x, G, I_1 \cup I_2) = N(x, G, I_1) + N(x, G, I_2) - N(x, G, I_1 \cap I_2).
\]

(3) Since $\mathbb{E}N\left(G,\rho,\displaystyle\bigcup_{i=1}^{\infty} I_i\right) \leq \mathbb{E}N(G,\rho,I) < +\infty$, we have:
\[
\mathbb{E}N\left(G,\rho,\bigcup_{i=1}^{\infty} I_i\right) = \int_{\R^n} \sum_{i=1}^{\infty} N(x, G, I_i) \rho(x)  dx = \sum_{i=1}^{\infty} \int_{\R^n} N(x, G, I_i) \rho(x)  dx
\]\[= \sum_{i=1}^{\infty} \mathbb{E}N(G,\rho,I_i).
\]

(4) and (5) are self-evident.

(6) If $I_1 \subseteq I_2$, then $N(x, G, I_1) \leq N(x, G, I_2)$, which implies:
\[
\Lambda_i(G, I_1) \subseteq \Lambda_i(G, I_2).
\]

(7) Since
\[
\mathbb{E}N(G,\rho,I) = \int_{\R^n} N(x, G, I) \rho(x)  dx = \int_{\R^n} \sum_{i=1}^{\infty} \mathbb{I}(x, \Lambda_i(G,I)) \rho(x)  dx < +\infty,
\]
we obtain:
\[
\mathbb{E}N(G,\rho,I) = \sum_{i=1}^\infty \int_{\R^n} \mathbb{I}(x, \Lambda_i(G,I)) \rho(x)  dx = \sum_{i=1}^\infty \int_{\Lambda_i(G,I)} \rho(x)  dx.
\]

(8) $\forall x\in \Xi(P(s,x),Q(s,x),J)$, there exists $\hat{s}\in J$ such that $P(\hat{s},x)=0$ and $Q(\hat{s},x)=0$. This implies:
\[
x\in L(\{P(s,x), Q(s,x)\},\{\hat{s}\})\subseteq L(\{P(s,x), Q(s,x)\}, J).
\]
\qed

We now define the cumulative expectation function:
\[
\varphi_{G,\rho}(t) = \mathbb{E}N(G,\rho,[l,t]), \quad t \in [l,u].
\]
Clearly, our objective is to compute $\varphi_{G,\rho}(u)$. By Property~\ref{pro-1}(4), $\varphi_{G,\rho}(t)$ is a monotonic increasing function. Additionally, we can establish its right-continuity.

\begin{lem}\label{lem-0}
    The function $\varphi_{G,\rho}(t)$ is right-continuous on $[l,u)$.
\end{lem}
\proof Let $t^* \in [l, u)$ be given. For any decreasing sequence $\{t_i\}_{i=1}^{+\infty} \subset (t^*, u]$ satisfying $\displaystyle\lim_{i \to +\infty} t_i = t^*$, we have:
\[
\varphi_{G,\rho}(t_1) - \varphi_{G,\rho}(t^*) = \mathbb{E}N(G, \rho, (t^*, t_1]) = \mathbb{E}N\left(G, \rho, \bigcup_{i=1}^{\infty} (t_{i+1}, t_i]\right)
\]
\[
= \sum_{i=1}^{\infty} \mathbb{E}N(G, \rho, (t_{i+1}, t_i]) = \sum_{i=1}^{\infty} (\varphi_{G,\rho}(t_i) - \varphi_{G,\rho}(t_{i+1})) = \varphi_{G,\rho}(t_1) - \lim_{i \to +\infty} \varphi_{G,\rho}(t_i).
\]
Therefore,
\[
\lim_{i \to \infty} \varphi_{G,\rho}(t_i) = \varphi_{G,\rho}(t^*),
\]
which establishes the right-continuity of $\varphi_{G,\rho}(t)$ on $[l, u)$. \qed

\begin{remark}
By Property~\ref{pro-1}(6), the function $\displaystyle\int_{\Lambda_i(G,[l,t])} \rho(x)  dx$ is monotonic increasing in $t$, and hence differentiable almost everywhere on $I$. Consequently,
\[
\varphi_{G,\rho}'(t) = \sum_{i=1}^{\infty} \left( \int_{\Lambda_i(G,[l,t])} \rho(x)  dx \right)'
\]
holds almost everywhere \cite{KolmogorovFomin}. Let $\Psi_{G,\rho}(I)$ denote the set of points where this equality fails; then $|\Psi_{G,\rho}(I)| = 0$.

It follows that
\[
\int_{I \setminus \Psi_{G,\rho}(I)} \varphi_{G,\rho}'(t)  dt = \int_{I} \varphi_{G,\rho}'(t)  dt \leq \varphi_{G,\rho}(u) - \varphi_{G,\rho}(l) = \mathbb{E}N(G,\rho,I) - \varphi_{G,\rho}(l).
\]
\end{remark}

The above analysis provides a lower bound for $\mathbb{E}N(G,\rho,I)$. However, a lower bound alone is insufficient; we seek to establish equality. This requires proving the absolute continuity of $\varphi_{G,\rho}(t)$ on $I$. Before addressing this, we introduce several lemmas necessary for the subsequent development.

\section{Lemmas Related to the Coarea Formula}
In this section, we present several lemmas required for the subsequent development. For brevity, we defer the proofs of these lemmas to the appendix.

\begin{lem}\cite{evans1992, federer1969}\label{coarea-1}
    Let $u(x) : \R^n \rightarrow \R^k$ ($k \leq n$) be a Lipschitz continuous function. Then for any integrable function $g(x) : \R^n \rightarrow \R$, the following holds:
    \[
    \int_{\R^n} g(x) \sqrt{\det(\nabla u(x)^\top \nabla u(x))}  dx = \int_{\R^k} dt \int_{u^{-1}(t)} g(x)  d\mathbb{H}^{n-k},
    \]
    where $\mathbb{H}^{n-k}$ denotes the $(n-k)$-dimensional Hausdorff measure, and $u^{-1}(t) = \{ x \in \R^n \mid u(x) = t \}$.
\end{lem}

Based on the above result, we derive the following corollary:

\begin{cor}\label{cor-1}
Let $P(s, x) : \hat{I} \times \R^n \rightarrow \R$, where $\hat{I} = [a, b]$.
Assume the following conditions hold:
\begin{itemize}
    \item[(1)] There exists a unique implicit function $s = u(x) : \Omega \rightarrow \R$ such that $P(s, x) = 0$ on $\hat{I} \times \Omega$, where $\Omega \subseteq L(P(s, x), \hat{I})$;
    \item[(2)] $P'(s, x)$ and $\nabla P(s, x)$ exist and are continuous on $\hat{I} \times \Omega$;
    \item[(3)] The function
          \[
          T(\hat{s}) = \int_{\Omega \cap L(P(s, x), \{\hat{s}\})} \rho(x) \frac{|P'(\hat{s}, x)|}{\|\nabla P(\hat{s}, x)\|}  d\mathbb{H}^{n-1}
          \]
          is Lebesgue integrable on $\hat{I}$;
    \item[(4)] $\displaystyle\int_{\Omega \cap \tau(\hat{I})} \rho(x)  dx = 0$, where $$\tau(\hat{I}) =\Xi(P(s, x), Q(s, x), \hat{I}), Q(s, x) = P'(s, x) \|\nabla P(s, x)\|;$$
    \item[(5)] For almost every $s \in \hat{I}$,
          \[
          \int_{\Omega \cap \tau(\{s\})} \rho(x) \frac{|P'(s, x)|}{\|\nabla P(s, x)\|}  d\mathbb{H}^{n-1} = 0;
          \]
    \item[(6)] $u(\Omega \setminus \tau(\hat{I})) = \hat{I}$.
\end{itemize}
Then the following hold:
\begin{itemize}
    \item[(1)] $\tau(\hat{I})=L(\{P(s,x),Q(s,x)\},\hat{I})=\Xi(P(s,x),Q(s,x),\hat{I})$;
    \item[(2)] \[
          \int_{\Omega} \rho(x)  dx = \int_{\hat{I}} d\hat{s} \int_{\Omega \cap L(P(s, x), \{\hat{s}\})} \rho(x) \frac{|P'(\hat{s}, x)|}{\|\nabla P(\hat{s}, x)\|}  d\mathbb{H}^{n-1}.
          \]
\end{itemize}
\end{cor}

\begin{remark}
    To avoid discussing the case where $\|\nabla P(s,x)\|=0$, the set $\Omega \cap \tau(\{s\})$ in condition (5) can be decomposed as:
    \[
    \Omega \cap \left( (\tau(\{s\})\cap\tau_1(\{s\})) \cup (\tau(\{s\}) \setminus \tau_1(\{s\})) \right),
    \]
    where $\tau_1(\{s\}) = L(\{P(s, x), \|\nabla P(s, x)\|\}, \{s\})$.
    Condition (5) can then be replaced by the following two conditions:
    \begin{itemize}
        \item[(5a)] $\displaystyle\int_{\Omega \cap (\tau(\{s\}) \setminus \tau_1(\{s\}))} \rho(x) \frac{|P'(s, x)|}{\|\nabla P(s, x)\|}  d\mathbb{H}^{n-1} = 0$;
        \item[(5b)] The Hausdorff dimension of $\Omega \cap (\tau(\{s\})\cap\tau_1(\{s\}))$ is less than $n-1$.
    \end{itemize}
\end{remark}

\begin{cor}\label{cor-2}
Let $P(s, x) : \hat{I} \times \R^n \rightarrow \R$, where $\hat{I} = [a, b]$ with $a\neq b$. Assume the following conditions hold:
\begin{itemize}
    \item[(1)] $P'(s, x)$ and $\nabla P(s, x)$ exist and are continuous on $\hat{I} \times L(P(s, x), \hat{I})$;
    \item[(2)] The function
          \[
          T(\hat{s}) = \int_{L(P(s, x), \{\hat{s}\})} \rho(x) \frac{|P'(\hat{s}, x)|}{\|\nabla P(\hat{s}, x)\|}  d\mathbb{H}^{n-1}
          \]
          is integrable on $\hat{I}$;
    \item[(3)] $\displaystyle\int_{\Xi(P(s,x), Q(s,x),  \hat{I})} \rho(x)  dx = 0$, where $Q(s, x) = P'(s, x) \|\nabla P(s, x)\|$.
\end{itemize}
Then
\[
\int_{L(P(s, x), \hat{I})} \rho(x)  dx \leq \int_{\hat{I}} d\hat{s} \int_{L(P(s, x), \{\hat{s}\})} \rho(x) \frac{|P'(\hat{s}, x)|}{\|\nabla P(\hat{s}, x)\|}  d\mathbb{H}^{n-1}.
\]
\end{cor}

\begin{cor}\label{cor-3}
Let $P(s, x) = \{P_i(s_i, x)\} : \hat{I}^2 \times \R^n \rightarrow \R^2$ for $i = 1, 2$, where $n \geq 2$ and $\hat{I} = [a, b]$ with $a\neq b$. Assume the following conditions hold:
\begin{itemize}
    \item[(1)] $P_i'(s_i, x)$ and $\nabla P_i(s_i, x)$ ($i = 1, 2$) exist and are continuous on $\hat{I} \times L(P(s, x), \hat{I}^2)$;
    \item[(2)] The function
          \[
          T(\hat{s}) = \int_{L(P(s, x), \{\hat{s}\})} \rho(x) \frac{|P_1'(\hat{s}_1, x)| |P_2'(\hat{s}_2, x)|}{\sqrt{V(\hat{s}, x)}}  d\mathbb{H}^{n-2}
          \]
          is integrable on $\hat{I}^2$, where
          \[
          V(s, x) = \|\nabla P_1(s_1, x)\|^2 \|\nabla P_2(s_2, x)\|^2 - (\nabla P_1(s_1, x)^\top \nabla P_2(s_2, x))^2;
          \]
    \item[(3)] $\displaystyle\int_{\Xi(P(s, x), Q(s, x), \hat{I}^2)} \rho(x)  dx = 0$, where $Q(s, x) = P'_1(s_1, x) P'_2(s_2, x) V(s, x)$.
\end{itemize}
Then
\[
\int_{L(P(s, x), \hat{I}^2)} \rho(x)  dx \leq \int_{\hat{I}^2} d\hat{s} \int_{L(P(s, x), \{\hat{s}\})} \rho(x) \frac{|P_1'(\hat{s}_1, x)| |P_2'(\hat{t}_2, x)|}{\sqrt{V(\hat{s}, x)}}  d\mathbb{H}^{n-2}.
\]
\end{cor}

\begin{lem}\cite{simon1983}\label{coarea-2}
Let $E \subseteq \R^n$ be an $m$-rectifiable Hausdorff measurable set, and let $u(x) : \R^n \rightarrow \R$ be a Lipschitz continuous function. Then for any nonnegative $\mathbb{H}^{m}$-measurable function $g : E \rightarrow \R$, the following holds:
\[
\int_{E} \|\nabla^E u(x)\| g(x)  d\mathbb{H}^{m} = \int_{-\infty}^{+\infty} dt \int_{u^{-1}(t) \cap E} g(x)  d\mathbb{H}^{m-1},
\]
where $\nabla^E u(x) = \displaystyle\sum_{i=1}^m \frac{\partial u(x)}{\partial W_i} W_i$, and $W_i \in \R^{n \times m}$ form an orthonormal basis for the tangent plane of $E$ at $x$.
\end{lem}

\begin{remark}
\begin{itemize}
    \item[(1)] Suppose $$E = \{ x \in \R^n \mid P_i(x) = 0,\ i = 1, \cdots, m \}$$ is an $m$-rectifiable Hausdorff measurable set with $\rank(\nabla P(x)) = m$. Let $W \in \R^{n \times (n-m)}$ satisfy $W^\top W = I_{n-m}$ and $\nabla P(x)^\top W = 0_{m \times (n-m)}$. Then
          \[
          \nabla^E u(x) = \sum_{i=1}^m \frac{\partial u(x)}{\partial W_i} W_i = W W^\top \nabla u(x),
          \]
          \[
          \|\nabla^E u(x)\|^2 = \nabla u(x)^\top W W^\top \nabla u(x).
          \]
          Let $\nabla u(x) = W \lambda + \nabla P(x) t$, where $\lambda \in \R^{n-m},\ t \in \R^m$. Then
          \[
          \|\nabla^E u(x)\|^2 = (t^\top \nabla P(x)^\top + \lambda^\top W^\top) W W^\top (W \lambda + \nabla P(x) t) = \|\lambda\|^2,
          \]
          \[
          \nabla P(x)^\top \nabla u(x) = \nabla P(x)^\top (W \lambda + \nabla P(x) t) = \nabla P(x)^\top \nabla P(x) t
          \]
          \[
          \Rightarrow t = (\nabla P(x)^\top \nabla P(x))^{-1} \nabla P(x)^\top \nabla u(x),
          \]
          \[
          \|\nabla u(x)\|^2 = \|\lambda\|^2 + t^\top \nabla P(x)^\top \nabla P(x) t.
          \]
          Therefore, we obtain:
          \[
          \|\nabla^E u(x)\| = \sqrt{ \|\nabla u(x)\|^2 - \nabla u(x)^\top \nabla P(x) (\nabla P(x)^\top \nabla P(x))^{-1} \nabla P(x)^\top \nabla u(x) }.
          \]

    \item[(2)] If $m = 1$, then
          \[
          \|\nabla^E u(x)\| = \frac{ \sqrt{ \|\nabla u(x)\|^2 \|\nabla P(x)\|^2 - (\nabla P(x)^\top \nabla u(x))^2 } }{\|\nabla P(x)\|}.
          \]
\end{itemize}
\end{remark}

\begin{cor}\label{cor-4}
Let $P_1(x) : \R^n \rightarrow \R$, $P_2(s, x) : \hat{I} \times \R^n \rightarrow \R$, where $\hat{I} = [a, b]$ with $a\neq b$ and $n \geq 2$. Assume the following conditions hold:
\begin{itemize}
    \item[(1)] $P_2'(s, x)$ and $\nabla P_2(s, x)$ exist and are continuous on $\hat{I} \times L(P_2(s, x), \hat{I})$, and $\nabla P_1(x)$ exists and is continuous on $L(P_2(s, x), \hat{I})$;
    \item[(2)] $\displaystyle\int_{\Omega \cap \Xi(P_2(s, x), Q(s, x), \hat{I})} \rho(x)  d\mathbb{H}^{n-1} = 0$, where
          \[
          \Omega = \{ x \in \R^n \mid P_1(x) = 0 \}, \quad Q(s, x) = P_2'(s, x) \|\nabla P_1(x)\| V(s, x),
          \]
          \[
          V(s, x) = \|\nabla P_2(s, x)\|^2 \|\nabla P_1(x)\|^2 - (\nabla P_1(x)^\top \nabla P_2(s, x))^2;
          \]
    \item[(3)] $\Omega$ is an $(n-1)$-rectifiable Hausdorff measurable set;
    \item[(4)] The function
          \[
          T(\hat{s}) = \int_{\Omega \cap L(P_2(s, x), \{\hat{s}\})} \rho(x) \frac{|P_2'(\hat{s}, x)| \|\nabla P_1(x)\|}{\sqrt{V(\hat{s}, x)}}  d\mathbb{H}^{n-2}
          \]
          is integrable on $\hat{I}$.
\end{itemize}
Then
\[
\int_{\Omega \cap L(P_2(s, x), \hat{I})} \rho(x)  d\mathbb{H}^{n-1} \leq \int_{\hat{I}} d\hat{s} \int_{\Omega \cap L(P_2(s, x), \{\hat{s}\})} \rho(x) \frac{|P_2'(\hat{s}, x)| \|\nabla P_1(x)\|}{\sqrt{V(\hat{s}, x)}}  d\mathbb{H}^{n-2}.
\]
\end{cor}

\begin{remark}
In the above corollaries, the set $\Xi(P(s,x),Q(s,x),J)$ is used extensively. However, $\Xi(P(s,x),Q(s,x),J)$ is generally more difficult to characterize than $L(\{P(s,x),Q(s,x)\},J)$. By Property \ref{pro-1}(8), we can control $\Xi(P(s,x),Q(s,x),J)$ through $L(\{P(s,x),Q(s,x)\},J)$. In many cases, these two sets are equal; however, for certain problems, $L(\{P(s,x),Q(s,x)\},J)$ might lead to an overestimation. For example:
\end{remark}

\begin{exam}
    Let $P(s,x)=sx_1+s^2x_2$, $Q(s,x)=(x_1+2sx_2)|s|\sqrt{1+s^2}$, $J=[0,1]$,
    $L(\{P(s,x),Q(s,x)\},J)=\R^2$, $\Xi(P(s,x),Q(s,x),J)=L(\{x_1+sx_2,x_1+2sx_2\},J)=\{0\}\times \R$. Clearly, the latter is much smaller than the former.
\end{exam}

\section{Absolute Continuity of $\varphi_{G,\rho}(t)$}
We now return to proving the absolute continuity of $\varphi_{G,\rho}(t)$.
By the Banach--Zaretsky theorem and the monotonicity of $\varphi_{G,\rho}(t)$, it suffices to prove that $\varphi_{G,\rho}(t)$ is continuous on $I$ and satisfies Luzin's N-property.

\begin{thm}\label{thm-1}
Assume the following conditions hold:
\begin{itemize}
    \item[(1)] $G'(t, x)$ and $\nabla G(t, x)$ exist and are continuous on $I \times L(G(t, x), I)$;
    \item[(2)] The function
          \begin{equation}\label{eq-T}
          T(\hat{t}) = \int_{L(G(t, x), \{\hat{t}\})} \rho(x) \frac{|G'(\hat{t}, x)|}{\|\nabla G(\hat{t}, x)\|}  d\mathbb{H}^{n-1}    
          \end{equation}
          is Lebesgue integrable on $I$;
    \item[(3)] For every closed subinterval $\bar{I}=[a,b]\subseteq I$ with $a\neq b$, 
          $$\displaystyle\int_{\Xi(G(t, x), Q_2(t, x), \bar{I})} \rho(x)  dx = 0,$$ 
          where $Q_2(t, x) = G'(t, x) \|\nabla G(t, x)\|$.
\end{itemize}
Then,  $\varphi_{G,\rho}(t)$ is absolutely continuous on $[l, u]$, and the following equality holds:
\[
\int_{I} \varphi_{G,\rho}'(t)  dt = \varphi_{G,\rho}(u) - \varphi_{G,\rho}(l) = \mathbb{E}N(G, \rho, I).
\]
\end{thm}

\proof
The proof is divided into several parts:

(1) First, we prove that for any subset $\hat{I} \subseteq I$ with $|\hat{I}| = 0$ (Lebesgue measure zero), we have:
\[
\int_{L(G(t, x), \hat{I})} \rho(x)  dx = 0.
\]
Since $|\hat{I}| = 0$, for each $k\in \mathbb{N}$, there exists a countable collection of open intervals $\{\tilde{I}_i^k\}$ (if $u\in \hat{I}$ or $l\in \hat{I}$, then at most two non-open intervals are needed) such that 
$$\hat{I} \subseteq \displaystyle \bigcup_i \tilde{I}^k_i \subseteq I,\quad \tilde{I}^k_i \cap \tilde{I}^k_j = \emptyset~ \text{for } i \neq j,\quad \left|\displaystyle\bigcup_i \tilde{I}^k_i\right| \leq \frac{1}{k},\quad \displaystyle\bigcup_i \tilde{I}^k_i\subseteq \displaystyle\bigcup_i \tilde{I}^{k-1}_i.$$

By Corollary~\ref{cor-2}, we have:
\[
\int_{L(G(t, x), \tilde{I}^k_i)} \rho(x)  dx \leq \int_{L(G(t, x), \mathrm{cl}(\tilde{I}^k_i))} \rho(x)  dx \leq \int_{\mathrm{cl}(\tilde{I}^k_i)} d\hat{t} \int_{L(G(t, x), \{\hat{t}\})} \rho(x) \frac{|G'(\hat{t}, x)|}{\|\nabla G(\hat{t}, x)\|}  d\mathbb{H}^{n-1}.
\]
Therefore,
\[
\int_{L(G(t, x), \hat{I})} \rho(x)  dx \leq  \int_{L(G(t, x), \cup_i \tilde{I}^k_i)} \rho(x)  dx \leq \sum_i  \int_{L(G(t, x), \tilde{I}^k_i)} \rho(x)  dx
\]
\[
\leq  \sum_i \int_{\mathrm{cl}(\tilde{I}^k_i)} d\hat{t} \int_{L(G(t, x), \{\hat{t}\})} \rho(x) \frac{|G'(\hat{t}, x)|}{\|\nabla G(\hat{t}, x)\|}  d\mathbb{H}^{n-1}\]\[
= \int_{\bigcup_i \mathrm{cl}(\tilde{I}^k_i)} d\hat{t} \int_{L(G(t, x), \{\hat{t}\})} \rho(x) \frac{|G'(\hat{t}, x)|}{\|\nabla G(\hat{t}, x)\|}  d\mathbb{H}^{n-1}. 
\]
The last equality holds because the integration domain differs by at most countably many points.

Since
\[
\lim_{k\rightarrow +\infty} \bigcup_i\mathrm{cl}(\tilde{I}^k_i)=\bigcap_k\bigcup_i\mathrm{cl}(\tilde{I}^k_i)\supseteq \hat{I}, \quad \left|\bigcap_k\bigcup_i\mathrm{cl}(\tilde{I}^k_i)\right|=0,
\]
we obtain:
\[
\int_{L(G(t, x), \hat{I})} \rho(x)  dx = \lim_{k\rightarrow +\infty} \int_{L(G(t,x),\cup_i I^k_i)} \rho(x)  dx 
\]
\[
=\lim_{k\rightarrow +\infty}\int_{\bigcup_i \mathrm{cl}(\tilde{I}^k_i)} d\hat{t} \int_{L(G(t, x), \{\hat{t}\})} \rho(x) \frac{|G'(\hat{t}, x)|}{\|\nabla G(\hat{t}, x)\|}  d\mathbb{H}^{n-1}
=0.
\]
This immediately implies:
\[
\mathbb{E}N(G, \rho, \hat{I}) = \int_{L(G(t, x), \hat{I})} N(x, G, \hat{I}) \rho(x)  dx \leq \int_{L(G(t, x), \hat{I})} N(x, G, I) \rho(x)  dx = 0,
\]

(2) Next, we prove that $\varphi_{G,\rho}(t)$ is continuous on $[l, u]$. By Lemma~\ref{lem-0}, $\varphi_{G,\rho}(t)$ is right-continuous on $[l, u)$; it remains to prove left-continuity.

Let $t^* \in (l, u]$ be given. For any increasing sequence $\{t_i\}_{i=1}^{+\infty} \subset [l, t^*)$ satisfying $\displaystyle\lim_{i \to +\infty} t_i = t^*$, and since $\mathbb{E}N(G, \rho, \{t\}) = 0$ for any $t \in I$, we have:
\[
\varphi_{G,\rho}(t^*) - \varphi_{G,\rho}(t_1) = \mathbb{E}N(G, \rho, (t_1, t^*]) = \mathbb{E}N(G, \rho, [t_1, t^*)).
\]
Consequently,
\[
\varphi_{G,\rho}(t^*) - \varphi_{G,\rho}(t_1) = \mathbb{E}N\left(G, \rho, \bigcup_{i=1}^{\infty} [t_{i}, t_{i+1})\right) = \sum_{i=1}^{\infty} \mathbb{E}N(G, \rho, [t_{i}, t_{i+1}))
\]
\[
= \sum_{i=1}^{\infty} \mathbb{E}N(G, \rho, (t_i, t_{i+1}]) = \sum_{i=1}^{\infty} (\varphi_{G,\rho}(t_{i+1}) - \varphi_{G,\rho}(t_i)) = \lim_{i \to +\infty} \varphi_{G,\rho}(t_i) - \varphi_{G,\rho}(t_1),
\]
which yields:
\[
\lim_{i \to +\infty} \varphi_{G,\rho}(t_i) = \varphi_{G,\rho}(t^*).
\]
This establishes the left-continuity of $\varphi_{G,\rho}(t)$ on $(l, u]$.
Therefore, $\varphi_{G,\rho}(t)$ is continuous on $[l, u]$.

(3) Finally, we prove that $\varphi_{G,\rho}(t)$ satisfies Luzin's N-property on $I$, i.e.,
\[
|\varphi_{G,\rho}(\hat{I})| = 0, \quad \text{for any } \hat{I} \subseteq I \text{ satisfying } |\hat{I}| = 0.
\]

Using the same sequence of intervals $\{\tilde{I}^k_i\}$ constructed for $\hat{I}$ in part (1), we have:
\[
0 \leq |\varphi_{G,\rho}(\hat{I})| \leq \left| \varphi_{G,\rho}\left( \bigcup_{i=1}^\infty \tilde{I}^k_i \right) \right| = \sum_{i=1}^{\infty} |\varphi_{G,\rho}( \tilde{I}^k_i )| = \sum_{i=1}^{\infty} |\varphi_{G,\rho}( \mathrm{cl}(\tilde{I}^k_i) )|\]
\[ = \sum_{i=1}^{\infty} (\varphi_{G,\rho}(b^k_i) - \varphi_{G,\rho}(a^k_i))
= \sum_i^{\infty} \mathbb{E}N(G, \rho, (a_i^k, b^k_i]) = \sum_i^{\infty} \mathbb{E}N(G, \rho, (a_i^k, b^k_i))\]
\[ = \mathbb{E}N\left(G, \rho, \bigcup_{i=1}^\infty \tilde{I}^k_i \right)
= \int_{\R^n} N\left(x, G, \bigcup_{i=1}^\infty \tilde{I}^k_i \right) \rho(x)  dx 
\]
\[= \int_{L\left(G(t, x), \bigcup_{i=1}^\infty \tilde{I}^k_i \right)} N\left(x, G, \bigcup_{i=1}^\infty \tilde{I}^k_i \right) \rho(x)  dx
\leq 
\int_{L\left(G(t, x), \bigcup_{i=1}^\infty \tilde{I}^k_i \right)} N\left(x, G, I \right) \rho(x)  dx.
\]
Since $\mathbb{E}N(G,\rho,I)=\displaystyle\int_{\R^n} N\left(x, G, I \right) \rho(x)  dx<+\infty$ and from the conclusion of part (1), we conclude that $|\varphi_{G,\rho}(\hat{I})| = 0$, i.e., Luzin's N-property holds on $I$.

Therefore, by the Banach--Zaretsky theorem, $\varphi_{G,\rho}(t)$ is absolutely continuous on $I$, and
\[
\int_{I} \varphi'_{G,\rho}(t)  dt = \varphi_{G,\rho}(u) - \varphi_{G,\rho}(l) = \mathbb{E}N(G, \rho, I).
\]
\qed

Following the same proof method, we obtain the following corollary.
\begin{cor}
    Given $t^*\in I$, assume the following conditions hold:
\begin{itemize}
    \item[(1)] $G'(t, x)$ and $\nabla G(t, x)$ exist and are continuous on $I \times L(G(t, x), I)$;
    \item[(2)] The function
          \begin{equation}\label{eq-T}
          T(\hat{t}) = \int_{L(G(t, x), \{\hat{t}\})} \rho(x) \frac{|G'(\hat{t}, x)|}{\|\nabla G(\hat{t}, x)\|}  d\mathbb{H}^{n-1}    
          \end{equation}
          is Lebesgue integrable on $I$;
    \item[(3)] there exists $\varepsilon>0, \forall \hat{\varepsilon}\in (0,\varepsilon)$, $$\displaystyle\int_{\Xi(G(t, x), Q_2(t, x), [t^*, t^*+\hat{\varepsilon}])} \rho(x)  dx = 0,$$ 
          where $Q_2(t, x) = G'(t, x) \|\nabla G(t, x)\|$.
\end{itemize}
Then, $$\int_{L(G(t,x),\{t^*\})}\rho(x)dx=0.$$
\end{cor}
\proof $\forall \hat{\varepsilon}\in (0,\varepsilon)$, $$\int_{L(G(t,x),\{t^*\})}\rho(x)dx\leq \int_{L(G(t,x),[t^*,t^*+\hat{\varepsilon}])}\rho(x)dx.$$
From Corollary \ref{cor-2}, we have
$$\int_{L(G(t,x),\{t^*\})}\rho(x)dx\leq \int_{[t^*,t^*+\hat{\varepsilon}]}T(\hat{t})d\hat{t},$$
which means $\displaystyle\int_{L(G(t,x),\{t^*\})}\rho(x)dx=0.$\qed

\section{Main Results}
We now proceed to compute $\varphi'_{G,\rho}(t)$.
By definition,
\[
\varphi'_{G,\rho}(t) = \lim_{\varepsilon \to 0_+} \frac{\mathbb{E}N(G,\rho,[l,t+\varepsilon]) - \mathbb{E}N(G,\rho,[l,t])}{\varepsilon} = \lim_{\varepsilon \to 0_+} \frac{\mathbb{E}N(G,\rho,(t,t+\varepsilon])}{\varepsilon}.
\]
This expression does not directly facilitate the computation of $\varphi'_{G,\rho}(t)$; we need to find more suitable computational schemes.

\begin{thm}\label{thm-2}
Given $t^* \in [l,u) \setminus \Psi_{G,\rho}(I)$, suppose $n\geq 2$ and there exists $\bar{\varepsilon} \in (0, u - t^*]$, depending on $t^*$, such that for every $\varepsilon \in (0, \bar{\varepsilon})$ the following conditions hold:
\begin{itemize}
    \item[(1)] $G'(t,x), G''(t,x), \nabla G(t,x), \nabla G'(t,x)$ exist and are continuous on $[t^*, t^*+\varepsilon] \times L(H((t_1,t_2),x), [t^*,t^*+\varepsilon]^2)$, where $H((t_1,t_2),x) = \{G(t_1,x), G'(t_2,x)\}$;
    \item[(2)] $\displaystyle\int_{L(G(t, x), \{t^*\})} \rho(x)  dx = 0$;
    \item[(3)] There exists a nonnegative, $(t^*, \varepsilon)$-parameterized, Lebesgue integrable function $C_{t^*, \varepsilon}(t_1,t_2)$ defined on $[t^*,t^*+\varepsilon]^2$ such that for almost every $(\hat{t}_1, \hat{t}_2) \in [t^*,t^*+\varepsilon]^2$,
          \begin{equation}\label{eq-K}
          K(\hat{t}_1,\hat{t}_2)=
          \int_{L(H((t_1,t_2),x), \{(\hat{t}_1,\hat{t}_2)\})} \rho(x) \frac{|G'(\hat{t}_1,x)| |G''(\hat{t}_2,x)|}{\sqrt{V((\hat{t}_1,\hat{t}_2),x)}}  d\mathbb{H}^{n-2} \leq C_{t^*, \varepsilon}(\hat{t}_1, \hat{t}_2),
          \end{equation}
          and
          \[
          \int_{[t^*,t^*+\varepsilon]} d\hat{t}_1 \int_{[t^*,t^*+\varepsilon]} C_{t^*, \varepsilon}(\hat{t}_1, \hat{t}_2)  d\hat{t}_2 = o(\varepsilon),
          \]
          where
          \[
          V((t_1,t_2),x) = \|\nabla G(t_1,x)\|^2 \|\nabla G'(t_2,x)\|^2 - (\nabla G(t_1,x)^\top \nabla G'(t_2,x))^2;
          \]
    \item[(4)] $\displaystyle\int_{\Xi(H((t_1,t_2),x), Q_1((t_1,t_2),x), [t^*,t^*+\varepsilon]^2)} \rho(x)  dx = 0$,
          where $$Q_1((t_1,t_2),x) = G'(t_1,x) G''(t_2,x) V((t_1,t_2),x).$$
\end{itemize}
Then,
\begin{equation}\label{eq-R}
\eta(t^*)=\left( \int_{\Theta(t^*,\varepsilon)} \rho(x)  dx \right)' \Big|_{\varepsilon=0}
\end{equation}
exists, and
\[
\varphi'_{G,\rho}(t^*) = \eta(t^*),
\]
where
\[
\Theta(t,\varepsilon) =\{ x \in \R^n \mid G(t,x)G(t+\varepsilon,x) \leq 0 \}.
\]
\end{thm}

\proof
By the definition,
\[
\varphi'_{G,\rho}(t^*) = \sum_{i=1}^{\infty} \left( \int_{\Lambda_i(G,(t^*,t^*+\varepsilon])} \rho(x)  dx \right)' \Big|_{\varepsilon=0},
\]
we need to calucale \[
\left( \int_{\Lambda_i(G,(t^*,t^*+\varepsilon])} \rho(x)  dx \right)' \Big|_{\varepsilon=0}, \forall i\in \N.
\]

We first prove that for $i \geq 2$,
\[
\left( \int_{\Lambda_i(G,(t^*,t^*+\varepsilon])} \rho(x)  dx \right)' \Big|_{\varepsilon=0} = 0.
\]

For any $x \in \Lambda_2(G,(t^*,t^*+\varepsilon])$, there exist $t_1, t_2 \in (t^*,t^*+\varepsilon]$ with $t_1 \neq t_2$ such that $G(t_1,x) = G(t_2,x) = 0$. This implies the existence of $t_3 \in (t^*,t^*+\varepsilon)$ such that $G'(t_3,x) = 0$. Therefore,
\[
\Lambda_2(G,(t^*,t^*+\varepsilon]) \subseteq L(G'(t,x),(t^*,t^*+\varepsilon)) \cap \Lambda_1(G,[t^*,t^*+\varepsilon]) \subseteq L(H((t_1,t_2),x),[t^*,t^*+\varepsilon]^2).
\]

By Corollary \ref{cor-3},
\[
\int_{L(H((t_1,t_2),x),[t^*,t^*+\varepsilon]^2)} \rho(x)  dx \leq \int_{[t^*,t^*+\varepsilon]} d\hat{t}_1 \int_{[t^*,t^*+\varepsilon]} K(\hat{t}_1,\hat{t}_2)d\hat{t}_2 
\]
\[
\leq \int_{[t^*,t^*+\varepsilon]} d\hat{t}_1 \int_{[t^*,t^*+\varepsilon]} C_{t^*, \varepsilon}(\hat{t}_1, \hat{t}_2)  d\hat{t}_2 = o(\varepsilon).
\]

Thus, by condition (2) we have:
\[
\int_{\Lambda_2(G,(t^*,t^*+\varepsilon])} \rho(x)  dx = o(\varepsilon),
\]
which implies
\[
\left( \int_{\Lambda_2(G,(t^*,t^*+\varepsilon])} \rho(x)  dx \right)' \Big|_{\varepsilon=0} = \left( \int_{L(H((t_1,t_2),x),[t^*,t^*+\varepsilon]^2)} \rho(x)  dx \right)' \Big|_{\varepsilon=0} = 0.
\]

Noting that
\[
\left( \int_{\Lambda_i(G,(t^*,t^*+\varepsilon])} \rho(x)  dx \right)' \Big|_{\varepsilon=0} \geq \left( \int_{\Lambda_{i+1}(G,(t^*,t^*+\varepsilon])} \rho(x)  dx \right)' \Big|_{\varepsilon=0},
\]
we conclude that
\[
\left( \int_{\Lambda_i(G,(t^*,t^*+\varepsilon])} \rho(x)  dx \right)' \Big|_{\varepsilon=0} = 0, \quad \text{for all } i \geq 2.
\]

We now consider
\[
\left( \int_{\Lambda_1(G,(t^*,t^*+\varepsilon])} \rho(x)  dx \right)' \Big|_{\varepsilon=0}.
\]

Observe that:
\[
\int_{\Lambda_1(G,(t^*,t^*+\varepsilon])} \rho(x)  dx = \int_{L(G'(t,x),(t^*,t^*+\varepsilon)) \cap \Lambda_1(G,(t^*,t^*+\varepsilon])} \rho(x)  dx\]
\[ + \int_{\Lambda_1(G,(t^*,t^*+\varepsilon]) \setminus L(G'(t,x),(t^*,t^*+\varepsilon))} \rho(x)  dx
= \int_{\Lambda_1(G,(t^*,t^*+\varepsilon]) \setminus L(G'(t,x),(t^*,t^*+\varepsilon))} \rho(x)  dx + o(\varepsilon).
\]
Therefore, the derivative
\[
\left( \int_{\Lambda_1(G,(t^*,t^*+\varepsilon]) \setminus L(G'(t,x),(t^*,t^*+\varepsilon))} \rho(x)  dx \right)' \Big|_{\varepsilon=0}
\]
exists, and
\[
\left( \int_{\Lambda_1(G,(t^*,t^*+\varepsilon])} \rho(x)  dx \right)' \Big|_{\varepsilon=0} = \left( \int_{\Lambda_1(G,(t^*,t^*+\varepsilon]) \setminus L(G'(t,x),(t^*,t^*+\varepsilon))} \rho(x)  dx \right)' \Big|_{\varepsilon=0}.
\]

By the condition (2),  we have:

\[
\left( \int_{\Lambda_1(G,(t^*,t^*+\varepsilon]) \setminus L(G'(t,x),(t^*,t^*+\varepsilon))} \rho(x)  dx \right)' \Big|_{\varepsilon=0}\]\[ = \left( \int_{\Lambda_1(G,[t^*,t^*+\varepsilon]) \setminus L(G'(t,x),(t^*,t^*+\varepsilon))} \rho(x)  dx \right)' \Big|_{\varepsilon=0},
\]
\[
=\left( \int_{\Lambda_1(G,(t^*,t^*+\varepsilon])} \rho(x)  dx \right)' \Big|_{\varepsilon=0} =\left( \int_{\Lambda_1(G,[t^*,t^*+\varepsilon])} \rho(x)  dx \right)' \Big|_{\varepsilon=0}.
\]
Note that
\[
\Theta(t^*,\varepsilon) \subseteq \Lambda_1(G,[t^*,t^*+\varepsilon]).
\]
If $\Lambda_1(G,[t^*,t^*+\varepsilon]) \setminus L(G'(t,x),(t^*,t^*+\varepsilon)) \neq \emptyset$, then for any $x$ in this set, we have $G'(t,x) \neq 0$ for all $t \in (t^*,t^*+\varepsilon)$, meaning $G(t,x)$ is strictly monotonic on $[t^*,t^*+\varepsilon]$. Consequently, $x \in \Theta(t^*,\varepsilon)$, and thus
\[
\Lambda_1(G,[t^*,t^*+\varepsilon]) \setminus L(G'(t,x),(t^*,t^*+\varepsilon)) \subseteq \Theta(t^*,\varepsilon) \subseteq \Lambda_1(G,[t^*,t^*+\varepsilon]).
\]
This inclusion also holds trivially if the set on the left is empty.

Therefore, 
\[\int_{\Lambda_1(G,[t^*,t^*+\varepsilon]) \setminus L(G'(t,x),(t^*,t^*+\varepsilon))} \rho(x)  dx \leq \int_{\Theta(t^*,\varepsilon)} \rho(x)  dx \leq \int_{\Lambda_1(G,[t^*,t^*+\varepsilon])} \rho(x)  dx
\]

Since
\[
\Lambda_1(G,[t^*,t^*]) \setminus L(G'(t,x),(t^*,t^*)) = \Theta(t^*,0) = \Lambda_1(G,[t^*,t^*]),
\]
we obtain:
\[
\left( \int_{\Lambda_1(G,[t^*,t^*+\varepsilon]) \setminus L(G'(t,x),(t^*,t^*+\varepsilon))} \rho(x)  dx \right)' \Big|_{\varepsilon=0} = \left( \int_{\Theta(t^*,\varepsilon)} \rho(x)  dx \right)' \Big|_{\varepsilon=0} \]\[= \left( \int_{\Lambda_1(G,[t^*,t^*+\varepsilon])} \rho(x)  dx \right)' \Big|_{\varepsilon=0}.
\]
Hence,
\[
\varphi'_{G,\rho}(t^*) = \left( \int_{\Theta(t^*,\varepsilon)} \rho(x)  dx \right)' \Big|_{\varepsilon=0}.
\]
\qed

\begin{remark}
\begin{itemize}
    \item[(1)] A typical choice for $C_{t^*, \varepsilon}(\hat{t}_1,\hat{t_2})$ is $C\varepsilon^{-\alpha}$, where $C>0$ and $\alpha<1$;
    \item[(2)] Theorem~\ref{thm-2} provides an expression for $\varphi'_{G,\rho}(t)$ under relatively weak conditions. However, it requires computing the derivative of a variable-domain integral, which is often not straightforward. To address this, we present Theorem~\ref{thm-3}. While requiring stronger conditions than Theorem~\ref{thm-2}, it avoids the differentiation of a variable-domain integral and provides a direct expression for $\varphi_{G,\rho}'(t)$.
\end{itemize}
\end{remark}
\begin{thm}\label{thm-3}
If all the conditions of Theorem~\ref{thm-2} except (2) are satisfied,  and the following additional conditions hold:
\begin{itemize}
    \item[(5)] $T(\hat{t})$ is Lebesgue integrable on $[t^*,t^*+\varepsilon]$;
    \item[(6)] $\displaystyle\int_{\Xi(G(t,x), Q_2(t,x), \tilde{I})} \rho(x)  dx = 0, \forall \tilde{I}=[t^*,\hat{\varepsilon}]\subseteq [t^*,t^*+\varepsilon]$ with $\hat{\varepsilon} \in (0,\bar{\varepsilon})$,
          where $$Q_2(t,x) = G'(t,x) \|\nabla G(t,x)\|;$$
    \item[(7)] For almost every $\hat{t} \in [t^*,t^*+\varepsilon]$,
          \[
          \int_{\Xi(G(t,x), Q_2(t,x), \{\hat{t}\})} \rho(x) \frac{|G'(\hat{t},x)|}{\|\nabla G(\hat{t},x)\|}  d\mathbb{H}^{n-1} = 0;
          \]
    \item[(8)] For almost every $\hat{t}_1 \in [t^*,t^*+\varepsilon]$,
          \[
          \int_{L(G(t,x), \{\hat{t}_1\}) \cap \Xi(G'(t_2,x), (Q_3)_{\hat{t}_1}(t_2,x), [t^*,t^*+\varepsilon])} \rho(x) \frac{|G'(\hat{t}_1,x)|}{\|\nabla G(\hat{t}_1,x)\|}  d\mathbb{H}^{n-1} = 0,
          \]
          where $(Q_3)_{\hat{t}_1}(t_2,x) = G''(t_2,x) \|\nabla G(\hat{t}_1,x)\| V((\hat{t}_1,t_2),x)$;
    \item[(9)] The following limits exist: $\exists \Omega_0\supseteq \{x\in \R^n| \rho(x)>0\}$ such that
          \[
          \lim_{\varepsilon \to 0_+} L(G'(t,x),(t^*,t^*+\varepsilon))\bigcap \Omega_0 = \Omega_1,
          \]
          \[
          \limsup_{\varepsilon \to 0_+} L(G(t,x),\{t^*+\varepsilon\})\bigcap \Omega_0 = \Omega_2,
          \]
          \[
          \limsup_{\varepsilon \to 0_+} \Xi(G(t,x), Q_2(s,x),  [t^*, t^*+\varepsilon])\bigcap \Omega_0 = \Omega_3,
          \]
          and satisfy:
          \[
          (L(G(t,x),\{t^*\})\bigcap \Omega_0) \setminus (\Omega_1 \cup \Omega_3) \neq \emptyset,\quad 
          \Omega_2 \setminus (\Omega_1 \cup \Omega_3) \neq \emptyset;
          \]
    \item[(10)] For almost every $\hat{t} \in [t^*,t^*+\varepsilon]$, $L(G(t,x),\{\hat{t}\})$ is an $(n-1)$-rectifiable Hausdorff measurable set;
    \item[(11)] $\beta = \displaystyle\lim_{\hat{t} \to t^*_+} T(\hat{t})$ exists.
\end{itemize}
Then
\[
\varphi'_{G,\rho}(t^*) = \beta.
\]
\end{thm}

\proof
We first prove that there exists $\hat{\varepsilon} \in (0, \varepsilon)$ such that for every $\tilde{\varepsilon} \in (0, \hat{\varepsilon}]$,
\[
(\Lambda_1(G,[t^*,t^*+\tilde{\varepsilon}])\bigcap \Omega_0) \setminus \left( L(G'(t,x),(t^*,t^*+\tilde{\varepsilon})) \cup \Xi(G(t,x), Q_2(t,x),[t^*,t^*+\tilde{\varepsilon}]) \right) \neq \emptyset.
\]

Suppose that it is not. Then there exists a sequence $\hat{\varepsilon}_k \to 0_+$ and corresponding $\tilde{\varepsilon}_k \in (0, \hat{\varepsilon}_k)$ such that
\[
L(G(t,x),[t^*,t^*+\tilde{\varepsilon}_k])\bigcap \Omega_0  \]\[\subseteq (L(G'(t,x),(t^*,t^*+\tilde{\varepsilon}_k)) \cup \Xi(G(t,x), Q_2(t,x),[t^*,t^*+\tilde{\varepsilon}_k]))\bigcap \Omega_0.
\]
Since
\[
\lim_{\varepsilon\rightarrow 0_+}L(G'(t,x),(t^*,t^*+\tilde{\varepsilon}_k))\bigcap \Omega_0= \Omega_1,
\]
\[
\lim_{\varepsilon\rightarrow 0_+}L(G(t,x),[t^*,t^*+\tilde{\varepsilon}_k])\bigcap \Omega_0= L(G(t,x),\{t^*\})\bigcap \Omega_0,
\]
\[
\limsup_{\varepsilon\rightarrow 0_+}\Xi(G(t,x), Q_2(t,x),[t^*,t^*+\tilde{\varepsilon}_k])\bigcap \Omega_0=\Omega_3,
\]
we would have
\[
L(G(t,x),\{t^*\})\bigcap \Omega_0 \subseteq \Omega_1 \cup \Omega_3,
\]
contradicting condition (9).

Therefore, there exists $\hat{\varepsilon} \in (0, \varepsilon)$ such that for every $\tilde{\varepsilon} \in (0, \hat{\varepsilon}]$,
\[
(\Lambda_1(G,[t^*,t^*+\tilde{\varepsilon}])\bigcap \Omega_0) \setminus \left( L(G'(t,x),(t^*,t^*+\tilde{\varepsilon})) \cup \Xi(G(t,x), Q_2(t,x),[t^*,t^*+\tilde{\varepsilon}]) \right) \neq \emptyset,
\]
which clearly implies
\[
(\Lambda_1(G,[t^*,t^*+\tilde{\varepsilon}]) \setminus L(G'(t,x),(t^*,t^*+\tilde{\varepsilon})))\bigcap \Omega_0 \neq \emptyset.
\]

For any $x \in \Lambda_1(G,[t^*,t^*+\hat{\varepsilon}]) \setminus L(G'(t,x),(t^*,t^*+\hat{\varepsilon}))$, we have $G'(t,x) \neq 0$ for all $t \in (t^*,t^*+\hat{\varepsilon})$, meaning $G(t,x)$ is strictly monotonic on $[t^*,t^*+\hat{\varepsilon}]$. Thus, $G(t,x) = 0$ has a unique solution $t$ in $[t^*,t^*+\hat{\varepsilon}]$, i.e., there exists a unique implicit function
\[
t = u(x) : (\Lambda_1(G,[t^*,t^*+\hat{\varepsilon}]) \setminus L(G'(t,x),(t^*,t^*+\hat{\varepsilon})))\bigcap \Omega_0 \to [t^*,t^*+\hat{\varepsilon}].
\]
This allows us to apply Corollary~\ref{cor-1} to compute
\[\left( \int_{(\Lambda_1(G,[t^*,t^*+\varepsilon]) \setminus L(G'(t,x),(t^*,t^*+\varepsilon)))\bigcap \Omega_0} \rho(x)  dx \right)' \Big|_{\varepsilon=0}\]\[=
\left( \int_{\Lambda_1(G,[t^*,t^*+\varepsilon]) \setminus L(G'(t,x),(t^*,t^*+\varepsilon))} \rho(x)  dx \right)' \Big|_{\varepsilon=0}.
\]

The only thing we need to prove is that there exists $\bar{\varepsilon} \in (0, \hat{\varepsilon})$ such that
\[
u\left( (\Lambda_1(G,[t^*,t^*+\bar{\varepsilon}])\bigcap \Omega_0) \setminus \left( L(G'(t,x),(t^*,t^*+\bar{\varepsilon})) \cup \Xi(G(t,x), Q_2(t,x),[t^*,t^*+\bar{\varepsilon}]) \right) \right)\] \[= [t^*,t^*+\bar{\varepsilon}].
\]

If for every $\bar{\varepsilon} \in (0, \hat{\varepsilon})$, there exists $\hat{t} \in [t^*,t^*+\bar{\varepsilon}]$ such that for all
\[
x \in (\Lambda_1(G,[t^*,t^*+\bar{\varepsilon}])\bigcap \Omega_0) \setminus \left( L(G'(t,x),(t^*,t^*+\bar{\varepsilon})) \cup \Xi(G(t,x), Q_2(t,x),[t^*,t^*+\bar{\varepsilon}]) \right)
\]
we have $G(\hat{t},x) \neq 0$ (i.e., $x \notin L(G(x,t),\{\hat{t}\})$), then
\[
(L(G(x,t),\{\hat{t}\})\bigcap \Omega_0) \subseteq (L(G'(t,x),(t^*,t^*+\bar{\varepsilon})) \cup \Xi(G(t,x), Q_2(t,x),[t^*,t^*+\bar{\varepsilon}]))\bigcap \Omega_0.
\]
Taking the limit as $\bar{\varepsilon} \to 0_+$, we get
\[
\Omega_2 \subseteq \Omega_1 \cup \Omega_3,
\]
contradicting condition (9).

Therefore,
\[
\int_{(\Lambda_1(G,[t^*,t^*+\bar{\varepsilon}])\bigcap \Omega_0) \setminus L(G'(t,x),(t^*,t^*+\bar{\varepsilon}))} \rho(x)  dx 
=\int_{\Lambda_1(G,[t^*,t^*+\bar{\varepsilon}]) \setminus L(G'(t,x),(t^*,t^*+\bar{\varepsilon}))} \rho(x)  dx \]\[= \int_{[t^*,t^*+\bar{\varepsilon}]} d\hat{t} \int_{L(G(t,x),\{\hat{t}\}) \setminus L(G'(t,x),(t^*,t^*+\bar{\varepsilon}))} \rho(x) \frac{|G'(\hat{t},x)|}{\|\nabla G(\hat{t},x)\|}  d\mathbb{H}^{n-1}.
\]
By condition (5) and (6), we obtain $\displaystyle\int_{L(G(t, x), \{t^*\})} \rho(x)  dx = 0$, so from Theorem~\ref{thm-2},
\[
\left( \int_{\Lambda_1(G,[t^*,t^*+\varepsilon]) \setminus L(G'(t,x),(t^*,t^*+\varepsilon))} \rho(x)  dx \right)' \Big|_{\varepsilon=0}
\]
exists, and we will prove that it is just $\beta.$

Note that
\[
\int_{[t^*,t^*+\varepsilon]} d\hat{t} \int_{L(G(t,x),\{\hat{t}\}) \setminus L(G'(t,x),(t^*,t^*+\varepsilon))} \rho(x) \frac{|G'(\hat{t},x)|}{\|\nabla G(\hat{t},x)\|}  d\mathbb{H}^{n-1}\]\[
= \int_{[t^*,t^*+\varepsilon]} T(\hat{t})  d\hat{t} - \int_{[t^*,t^*+\varepsilon]} d\hat{t} \int_{L(G(t,x),\{\hat{t}\}) \cap L(G'(t,x),(t^*,t^*+\varepsilon))} \rho(x) \frac{|G'(\hat{t},x)|}{\|\nabla G(\hat{t},x)\|}  d\mathbb{H}^{n-1}.
\]

By Corollary~\ref{cor-4},
\[
\int_{L(G(t,x),\{\hat{t}\}) \cap L(G'(t_2,x),(t^*,t^*+\varepsilon))} \rho(x) \frac{|G'(\hat{t},x)|}{\|\nabla G(\hat{t},x)\|}  d\mathbb{H}^{n-1}
\]\[\leq \int_{L(G(t,x),\{\hat{t}\}) \cap L(G'(t_2,x),[t^*,t^*+\varepsilon])} \rho(x) \frac{|G'(\hat{t},x)|}{\|\nabla G(\hat{t},x)\|}  d\mathbb{H}^{n-1}
\]
\[
\leq \int_{[t^*,t^*+\varepsilon]} d\hat{t}_2 \int_{L(G(t,x),\{\hat{t}\}) \cap L(G'(t_2,x),\{\hat{t}_2\})} \rho(x) \frac{|G'(\hat{t},x)| |G''(\hat{t}_2,x)|}{\sqrt{V((\hat{t},\hat{t}_2),x)}}  d\mathbb{H}^{n-2}\] \[= \int_{[t^*,t^*+\varepsilon]} K(\hat{t},\hat{t}_2)d\hat{t}_2.
\]

Thus,
\[
\int_{[t^*,t^*+\varepsilon]} d\hat{t} \int_{L(G(t,x),\{\hat{t}\}) \cap L(G'(t,x),(t^*,t^*+\varepsilon))} \rho(x) \frac{|G'(\hat{t},x)|}{\|\nabla G(\hat{t},x)\|}  d\mathbb{H}^{n-1} \]
\[\leq \int_{[t^*,t^*+\varepsilon]}d\hat{t}\int_{[t^*,t^*+\varepsilon]} K(\hat{t},\hat{t}_2)d\hat{t}_2 
\leq \int_{[t^*,t^*+\varepsilon]}d\hat{t}\int_{[t^*,t^*+\varepsilon]} C_{t^*, \varepsilon}(\hat{t},\hat{t}_2)d\hat{t}_2 = o(\varepsilon).
\]
Therefore,
\[
\left( \int_{[t^*,t^*+\varepsilon]} d\hat{t} \int_{L(G(t,x),\{\hat{t}\}) \cap L(G'(t,x),(t^*,t^*+\varepsilon))} \rho(x) \frac{|G'(\hat{t},x)|}{\|\nabla G(\hat{t},x)\|}  d\mathbb{H}^{n-1} \right)' \Big|_{\varepsilon=0} = 0.
\]
Hence,
\[
\left( \int_{\Lambda_1(G,[t^*,t^*+\varepsilon]) \setminus L(G'(t,x),(t^*,t^*+\varepsilon))} \rho(x)  dx \right)' \Big|_{\varepsilon=0} = \left( \int_{[t^*,t^*+\varepsilon]} T(\hat{t})  d\hat{t} \right)' \Big|_{\varepsilon=0}.
\]
Denote this common value by $\gamma$.

If $\beta \neq \gamma$, without loss of generality assume $\gamma > \beta$. Then there exists $\tilde{\varepsilon} > 0$ such that for $0 < \varepsilon < \tilde{\varepsilon}$,
\[
|T(t^*+\varepsilon) - \beta| < \frac{\gamma - \beta}{2}, \quad \left| \frac{\displaystyle\int_{[t^*,t^*+\varepsilon]} T(\hat{t})  d\hat{t}}{\varepsilon} - \gamma \right| < \frac{\gamma - \beta}{2},
\]
which implies
\[
\int_{[t^*,t^*+\tilde{\varepsilon}]} T(\hat{t})  d\hat{t} < \frac{\beta + \gamma}{2} \tilde{\varepsilon} < \int_{[t^*,t^*+\tilde{\varepsilon}]} T(\hat{t})  d\hat{t},
\]
a contradiction. Therefore,
\[
\left( \int_{[t^*,t^*+\varepsilon]} T(\hat{t})  d\hat{t} \right)' \Big|_{\varepsilon=0} = \beta.
\]

We conclude that
\[
\varphi'_{G,\rho}(t^*) = \beta.
\]
\qed

Combine all the above theorems, we could obtain the following corollaries.
\begin{cor}\label{cor-main1}
Under the conditions of Theorem~\ref{thm-2}, and with the following additional assumptions:
\begin{itemize}
    \item[(1)] $G'(t,x), G''(t,x), \nabla G(t,x), \nabla G'(t,x)$ exist and are continuous on $I \times L(H((t_1,t_2),x), I^2)$, where $H((t_1,t_2),x) = \{G(t_1,x), G'(t_2,x)\}$;
    \item[(2)] The function $T(\hat{t})$ is Lebesgue integrable on $I$;
    \item[(3)] For every $\hat{t} \in I$, there exists $\varepsilon > 0$ and an $(\hat{t}, \varepsilon)$-parameterized Lebesgue integrable function $C_{\hat{t}, \varepsilon}(\hat{t}_1, \hat{t}_2)$ defined on $[\hat{t},\hat{t}+\varepsilon]^2$ such that for almost every $(\hat{t}_1, \hat{t}_2) \in [\hat{t},\hat{t}+\varepsilon]^2$,
          \[
          K(\hat{t}_1,\hat{t}_2) \leq C_{\hat{t}, \varepsilon}(\hat{t}_1, \hat{t}_2),
          \]
          and
          \[
          \int_{[\hat{t},\hat{t}+\varepsilon]} d\hat{t}_1 \int_{[\hat{t},\hat{t}+\varepsilon]} C_{\hat{t}, \varepsilon}(\hat{t}_1, \hat{t}_2)  d\hat{t}_2 = o(\varepsilon);
          \]
    \item[(4)] For every closed subinterval $\bar{I}=[a,b] \subseteq I$ with $a \neq b$,
          \[
          \int_{\Xi(G(t,x), Q_2(t,x), \bar{I})} \rho(x)  dx = 0,
          \]
          where $Q_2(t,x) = G'(t,x) \|\nabla G(t,x)\|$;
    \item[(5)] For almost every $\hat{t} \in I$, there exists $\varepsilon > 0$ such that
          \[
          \int_{\Xi(H((t_1,t_2),x), Q_1((t_1,t_2),x), [\hat{t},\hat{t}+\varepsilon]^2)} \rho(x)  dx = 0,
          \]
          where $Q_1((t_1,t_2),x) = G'(t_1,x) G''(t_2,x) V((t_1,t_2),x)$.
\end{itemize}
Then
\[
\mathbb{E}N(G, \rho, I) = \int_{I} \eta(t)  dt,
\]
where $\eta(t)$ is defined in \eqref{eq-R}.
\end{cor}

\begin{cor}\label{cor-main2}
If conditions (1), (3)--(5) of Corollary~\ref{cor-main1} are satisfied, $n\geq 2$ and the following additional conditions hold:
\begin{itemize}
    \item[(2')] $T(\hat{t})$ is continuous almost everywhere on $I$;
    \item[(6)] For almost every $\hat{t} \in I$,
          \[
          \int_{\Xi(G(t,x), Q_2(t,x), \{\hat{t}\})} \rho(x) \frac{|G'(\hat{t},x)|}{\|\nabla G(\hat{t},x)\|}  d\mathbb{H}^{n-1} = 0;
          \]
    \item[(7)] For almost every $\hat{t} \in I$, there exists $\varepsilon > 0$ such that for almost every $\hat{t}_1\in [\hat{t},\hat{t}+\varepsilon]$,
          \[
          \int_{L(G(t,x), \{\hat{t}_1\}) \cap \Xi(G'(t_2,x), (Q_3)_{\hat{t}_1}(t_2,x), [\hat{t},\hat{t}+\varepsilon])} \rho(x) \frac{|G'(\hat{t}_1,x)|}{\|\nabla G(\hat{t}_1,x)\|}  d\mathbb{H}^{n-1} = 0,
          \]
          where $(Q_3)_{\hat{t}_1}(t_2,x) = G''(t_2,x) \|\nabla G(\hat{t}_1,x)\| V((\hat{t}_1,t_2),x)$;
    \item[(8)] For almost every $\hat{t} \in I$,
          \[
          \lim_{\varepsilon \to 0_+} L(G'(t,x),(\hat{t},\hat{t}+\varepsilon))\bigcap \Omega_0 = \Omega_1,
          \]
          \[
          \limsup_{\varepsilon \to 0_+} L(G(t,x),\{\hat{t}+\varepsilon\})\bigcap \Omega_0 = \Omega_2,
          \]
          \[
          \limsup_{\varepsilon \to 0_+} \Xi(G(t,x), Q_2(t,x), [\hat{t},\hat{t}+\varepsilon])\bigcap \Omega_0 = \Omega_3,
          \]
          and
          \[
          (L(G(t,x),\{\hat{t}\})\bigcap \Omega_0) \setminus (\Omega_1 \cup \Omega_3) \neq \emptyset, \quad
          \Omega_2 \setminus (\Omega_1 \cup \Omega_3) \neq \emptyset;
          \]
    \item[(9)] For almost every $\hat{t} \in I$, the set $L(G(t,x),\{\hat{t}\})$ is an $(n-1)$-rectifiable Hausdorff measurable set.
\end{itemize}
Then
\[
\mathbb{E}N(G, \rho, I) = \int_{I} T(t)  dt.
\]
\end{cor}

\begin{remark}\label{remark-conditions}
Many of the conditions in Corollaries~\ref{cor-main1} and~\ref{cor-main2} can be replaced by stronger but more easily verifiable conditions:
\begin{itemize}
    \item[(3')] For almost every $t \in I$, there exists $\varepsilon > 0$ and $C > 0$ such that for almost every $(\hat{t}_1, \hat{t}_2) \in [t,t+\varepsilon]^2$, we have $K(\hat{t}_1,\hat{t}_2) \leq C$;
    \item[(4')] $\displaystyle\int_{L(\{G(t,x), Q_2(t,x)\}, I)} \rho(x)  dx = 0$;
    \item[(4'')] The Hausdorff dimension of $L(\{G(t,x), Q_2(t,x)\}, I)$ is less than $n$;
    \item[(5')] For almost every $\hat{t} \in I$, there exists $\varepsilon > 0$ such that
          \[
          \int_{L(\{H((t_1,t_2),x), Q_1((t_1,t_2),x)\}, [\hat{t},\hat{t}+\varepsilon]^2)} \rho(x)  dx = 0;
          \]
    \item[(5'')] For almost every $\hat{t} \in I$, there exists $\varepsilon > 0$ such that the Hausdorff dimension of $L(\{H((t_1,t_2),x), Q_1((t_1,t_2),x)\}, [\hat{t},\hat{t}+\varepsilon]^2)$ is less than $n$;
    \item[(6')] For almost every $\hat{t} \in I$, the Hausdorff dimension of $L(\{G(t,x), Q_2(t,x)\}, \{\hat{t}\})$ is less than $n-1$;
    \item[(7')] For almost every $\hat{t} \in I$, there exists $\varepsilon > 0$ such that for almost every $\hat{t}_1 \in [\hat{t},\hat{t}+\varepsilon]$,
          \[
          \int_{L(G(t,x), \{\hat{t}_1\}) \cap L(\{G'(t_2,x), (Q_3)_{\hat{t}_1}(t_2,x)\}, [\hat{t},\hat{t}+\varepsilon])} \rho(x) \frac{|G'(\hat{t}_1,x)|}{\|\nabla G(\hat{t}_1,x)\|}  d\mathbb{H}^{n-1} = 0;
          \]
    \item[(7'')] For almost every $\hat{t} \in I$, there exists $\varepsilon > 0$ such that for almost every $\hat{t}_1 \in [\hat{t},\hat{t}+\varepsilon]$, the Hausdorff dimension of $$L(G(t,x),\{\hat{t}_1\}) \cap L(\{G'(t_2,x), (Q_3)_{\hat{t}}(t_2,x)\}, [\hat{t},\hat{t}+\varepsilon])$$ is less than $n-1$;
    \item[(8')] For almost every $\hat{t} \in I$,
          \[
          \lim_{\varepsilon \to 0_+} L(G'(t,x),(\hat{t},\hat{t}+\varepsilon)) = \Omega_1,
          \]
          \[
          \lim_{\varepsilon \to 0_+} L(G(t,x),\{\hat{t}+\varepsilon\}) = L(G(t,x),\{\hat{t}\}),
          \]
          and
          \[
          L(G(t,x),\{\hat{t}\}) \setminus (\Omega_1 \cup L(\{G(t,x), Q_2(t,x)\},\{\hat{t}\})) \neq \emptyset.
          \]
\end{itemize}
\end{remark}

\section{Linear Case: $G(t,x) = f(t)^\top x - f_0(t)$}
In this section, we analyze the important special case of linear functions. Due to the special structure of linear functions, we can derive more explicit expressions for $T(t)$.

\begin{lem}
Let $G(t,x) = f(t)^\top x - f_0(t)$, where $f(t) = (f_1(t), \dots, f_n(t))^\top$, and let $t^* \in I$. If $f(t)$ and $f_0(t)$ are continuously differentiable in a neighborhood of $t^*$, and the matrix $D=\begin{pmatrix} f(t^*)^\top \\ f'(t^*)^\top \end{pmatrix}$ has full rank, then
\[
T(t^*) = \frac{S^{3/2}}{\|f(t^*)\|} \int_{\R^{n-1}} \rho(W + Z\mu) |\mu_1|  d\mu,
\]
where $S$, $W$, and $Z$ are defined as follows:
\[
W = D^\top(DD^\top)^{-1}\begin{pmatrix} f_0(t^*) \\ f'_0(t^*) \end{pmatrix}, \quad
Z = \begin{pmatrix} f'(t^*) - \displaystyle\frac{f(t^*)^\top f'(t^*)}{f(t^*)^\top f(t^*)} f(t^*) & J \end{pmatrix},
\]
and
\[
S = \frac{f'(t^*)^\top f'(t^*) f(t^*)^\top f(t^*) - (f'(t^*)^\top f(t^*))^2}{f(t^*)^\top f(t^*)},
\]
where $J =Null(D)\in \R^{n \times (n-2)}$ be such that $D J = 0_{2 \times (n-2)}$ and $J^\top J = I_{(n-2) \times (n-2)}$. 

Furthermore, if $\rho(x) = \hat{\rho}(x^\top x)$ is a radial function, then
\[
T(t^*) = \frac{S^{3/2}}{\|f(t^*)\|} \int_{\R^{n-1}} \hat{\rho}(U + 2q\mu_1 + S\mu_1^2 + \|\mu_{2:n-1}\|^2) |\mu_1|  d\mu,
\]
where $U$ and $q$ are defined by:
\[
U = W^\top W = \frac{\|f_0(t^*) f'(t^*) - f'_0(t^*) f(t^*)\|^2}{S \|f(t^*)\|^2},
\]
\[
q = f'_0(t^*) - \frac{f(t^*)^\top f'(t^*)}{f(t^*)^\top f(t^*)} f_0(t^*).
\]
\end{lem}

\proof
By definition,
\[
T(t^*) = \int_{L(G(t,x),\{t^*\})} \frac{\rho(x) |G'(t^*,x)|}{\|\nabla_x G(t^*,x)\|}  d\mathbb{H}^{n-1} 
= \int_{L(G(t,x),\{t^*\})} \frac{\rho(x) |f'(t^*)^\top x - f'_0(t^*)|}{\|f(t^*)\|}  d\mathbb{H}^{n-1}.
\]

Then the level set can be parameterized as
\[
L(G(t,x),\{t^*\}) = \{ x \in \R^n \mid x = W + Z\mu, \mu \in \R^{n-1} \}.
\]

Direct computation shows that
\[
Z^\top Z = \begin{pmatrix} S & 0_{1 \times (n-2)} \\ 0_{(n-2) \times 1} & I_{(n-2) \times (n-2)} \end{pmatrix},
\]
and
\[
f'(t^*)^\top Z = (S, 0, \dots, 0), \quad f(t^*)^\top W = f_0(t^*).
\]

The Jacobian factor is
\[
\frac{|G'(t^*,x)|}{\|\nabla G(t^*,x)\|} = \frac{|f'(t^*)^\top x - f'_0(t^*)|}{\|f(t^*)\|} = \frac{|f'(t^*)^\top (W + Z\mu) - f'_0(t^*)|}{\|f(t^*)\|} = \frac{|S\mu_1|}{\|f(t^*)\|}.
\]

Therefore,
\[
T(t^*) = \frac{1}{\|f(t^*)\|} \int_{\R^{n-1}} \rho(W + Z\mu) |S\mu_1| \sqrt{\det(Z^\top Z)}  d\mu\]\[
= \frac{S^{3/2}}{\|f(t^*)\|} \int_{\R^{n-1}} \rho(W + Z\mu) |\mu_1|  d\mu.
\]

If $\rho(x) = \hat{\rho}(x^\top x)$ is radial, then
\[
x^\top x = (W + Z\mu)^\top (W + Z\mu) = W^\top W + 2W^\top Z\mu + \mu^\top Z^\top Z\mu.
\]
Computing $W^\top Z$ and $\mu^\top Z^\top Z\mu$ yields the expression involving $U$, $q$, and $S$.
\qed

\begin{remark}
We now consider several specific examples:
\begin{itemize}
    \item[(1)] If $x \sim N(0_n, I_{n\times n})$ follows a standard normal distribution, then
          \[
          T(t^*) = \frac{S^{3/2}}{\sqrt{(2\pi)^n} \|f(t^*)\|} \int_{\R^{n-1}} \exp\left(-\frac{U + 2q\mu_1 + S\mu_1^2 + \|\mu_{2:n-1}\|^2}{2}\right) |\mu_1|  d\mu
          \]
          \[
          = \frac{S^{3/2}}{(2\pi)^{n/2} \|f(t^*)\|} \int_{\R} \exp\left(-\frac{U + 2q\mu_1 + S\mu_1^2}{2}\right) |\mu_1|  d\mu_1 \cdot (2\pi)^{(n-1)/2}
          \]
          \[
          = \frac{\sqrt{S}}{(2\pi)^{1/2} \|f(t^*)\|} \int_{\R} \exp\left(-\frac{U + 2q\mu_1 + S\mu_1^2}{2}\right) |\mu_1|  d\mu_1.
          \]
          This integral can be evaluated in closed form.

    \item[(2)] If $x \sim U(B(0_n,r))$ is uniformly distributed on the ball of radius $r$, then
          \[
          T(t^*) = \frac{S^{3/2}}{|B(0_n, r)| \|f(t^*)\|} \int_{\Psi} |\mu_1|  d\mu,
          \]
          where $\Psi = \{\mu \in \R^{n-1} \mid U + 2q\mu_1 + S\mu_1^2 + \|\mu_{2:n-1}\|^2 \leq r^2\}$.
          This can be computed as
          \[
          T(t^*) = \frac{S^{3/2}}{|B(0_n, r)| \|f(t^*)\|} \int_{p_1}^{p_2} |\mu_1|  d\mu_1 \int_{\hat{\Psi}(\mu_1)} d\mu_{2:n-1},
          \]
          where $\hat{\Psi}(\mu_1) = \{\mu_{2:n-1} \in \R^{n-2} \mid \|\mu_{2:n-1}\|^2 \leq r^2 - (U + 2q\mu_1 + S\mu_1^2)\}$,
          and $p_1, p_2$ are the roots of $U + 2q\mu_1 + S\mu_1^2 = r^2$.

    \item[(3)] A special case occurs when $q = 0$, which is equivalent to $f_0(t) \equiv C \|f(t)\|$ for some constant $C$. In this case, the expressions simplify significantly.
\end{itemize}
\end{remark}

\begin{exam}
Let $I = [0, 1]$, $f_i(t) = t^{i-1}$ for $i = 1, \dots, n$, $f_0(t) = 0$, and $x \sim N(0_n, I_{n\times n})$. Then
\[
T(t) = \frac{\sqrt{1 - n^2 t^{2n-2} + t^{4n} + (2n^2 - 2)t^{2n} - n^2 t^{2n+2}}}{\pi (t^{2n} - 1)(t^2 - 1)},
\]
and
\[
\mathbb{E}N(G,\rho,I) = \int_{0}^{1} \frac{\sqrt{1 - n^2 t^{2n-2} + t^{4n} + (2n^2 - 2)t^{2n} - n^2 t^{2n+2}}}{\pi (t^{2n} - 1)(t^2 - 1)}  dt.
\]
\end{exam}

\begin{exam}
Let $I = [0, 1]$, $f_i(t) = t^{i-1}$ for $i = 1, \dots, n$, $f_0(t) = 0$, and $x \sim U(B(0_n,r))$. Then $T(t)$ and $\mathbb{E}N(G,\rho,I)$ are the same as in the previous example.
\end{exam}

\begin{exam}
Let $I = [0, +\infty)$, $f_i(t) = \sqrt{\binom{n}{i-1}} t^{i-1}$ for $i = 1, \dots, n+1$, $f_0(t) = \mu$, and $x \sim N(0_n, I_{n\times n})$. Then
\[
T(t)=\frac{\sqrt{n}}{\pi(1+t^2)}\exp\left(-\frac{\mu^2 (1+nt^2)}{2(1+t^2)^n}\right)\left(1+\frac{\sqrt{2\pi}}{2}\kappa(t)\left(2\Phi(\kappa(t))-1\right)\exp(\kappa^2(t))\right),
\]
where $\kappa(t)=\displaystyle\frac{-\mu \sqrt{n}t}{(1+t^2)^{n/2}}$, and $\Phi(t)$ is the cumulative distribution function of the standard normal distribution. In particular, when $\mu=0$,
\[
T(t) = \frac{\sqrt{n}}{\pi(1+t^2)}, \quad
\mathbb{E}N(G,\rho,I) = \int_{0}^{+\infty} \frac{\sqrt{n}}{\pi(1+t^2)}  dt=\frac{\sqrt{n}}{2}.
\]
\end{exam}

\begin{exam}
Let $I = [0, +\infty)$, $f_i(t) = \sqrt{\binom{n}{i-1}} t^{i-1}$ for $i = 1, \dots, n+1$, $f_0(t) = (1+t^2)^{n/2}$, and $x \sim U(B(0_n,r))$. Then
\[
T(t)= \begin{cases} 
      \displaystyle\frac{\sqrt{n}}{\pi(1+t^2)}\left(1-\frac{1}{r^2(1+t^2)^n}\right)^{n/2}, & t \geq \sqrt{\max\{r^{-2/n}-1, 0\}} \\
      0, & \text{otherwise}
  \end{cases}.
\]
In particular, when $r\geq 1$ and $n=2k$,
\[
\mathbb{E}N(G,\rho,I) = \sqrt{\frac{k}{2\pi}}\sum_{i=0}^k\frac{(-1)^i\binom{k}{i}}{r^{2i}}\frac{\Gamma(2ki+\frac{1}{2})}{\Gamma(2ki+1)}.
\]
\end{exam}

\begin{exam}
Let $I = [0, 2\pi]$, $f_{2i-1}(t) = \cos(it)$, $f_{2i}(t) = \sin(it)$ for $i = 1, \dots, n$, $f_0(t) = \mu$, and $x \sim N(0_n, I_{n\times n})$. Then
\[
T(t) = \sqrt{\frac{(n+1)(2n+1)}{6}} \exp\left(-\frac{\mu^2}{2}\right),
\]
\[
\mathbb{E}N(G,\rho,I) = 2\pi \sqrt{\frac{(n+1)(2n+1)}{6}} \exp\left(-\frac{\mu^2}{2}\right).
\]
\end{exam}

\begin{exam}
Let $I = [0, 2\pi]$, $f_{2i-1}(t) = \cos(it)$, $f_{2i}(t) = \sin(it)$ for $i = 1, \dots, n$, $f_0(t) = \mu$, and $x \sim U(B(0_n,r))$. Then
\[
T(t)= \begin{cases}
      \displaystyle\sqrt{\frac{(n+1)(2n+1)}{6}} \left(1 - \frac{\mu^2}{r^2}\right)^n, & |\mu| \leq r \\
      0, & \text{otherwise}
  \end{cases},
\]
\[
\mathbb{E}N(G,\rho,I) = \begin{cases}
      2\pi \displaystyle\sqrt{\frac{(n+1)(2n+1)}{6}} \left(1 - \frac{\mu^2}{r^2}\right)^n, & |\mu| \leq r \\
      0, & \text{otherwise}
  \end{cases}.
\]
\end{exam}

We conclude this section with two simple nonlinear examples that can be treated using our method.

\begin{exam}
Let $I=[l,u]$ with $l\geq 1$, $n\geq 3$, $G(t,x) = x_1^2 - \displaystyle\sum_{i=2}^n t^{i-1} x_i^2$, and $x_i \sim U(-1,1)$. Then
\[
T(t)=\frac{C_n}{2^n n t^{\frac{n(n-1)}{4}+1}}, \quad
\mathbb{E}N(G,\rho,I)=\frac{4C_n}{2^n n^2 (n-1)}\left(l^{-\frac{n(n-1)}{4}} - u^{-\frac{n(n-1)}{4}}\right),
\]
where $C_n$ is a constant given by
\[
C_3=\pi, \quad C_n=\int_{[0,\pi]^{n-3}}\left(\prod_{i=1}^{n-3}\sin^{n-2-i}(z_i)\right)\left(\sum_{i=1}^{n-2}a_i \prod_{j=1}^{n-2-i}\sin^2(z_j)\right)dz,
\]
with $a_1=3\pi$, $a_i=2\pi$ for $i=2,\dots,n-2$, and $n\geq 4$. For example, $C_4=8\pi$, $C_5=5\pi^2$, $C_6=8\pi^2$.
\end{exam}

\begin{exam}
Let $I=[l,u]$ with $l\geq 1$, $G(t,x)=\displaystyle\sum_{i=1}^n t^{i-1} x_i^2 - 1$, and $x \sim U(B(0_n,1))$. Then
\[
T(t)=\frac{n(n-1)}{4}t^{-\frac{n(n-1)}{4}-1}, \quad
\mathbb{E}N(G,\rho,I)=l^{-\frac{n(n-1)}{4}} - u^{-\frac{n(n-1)}{4}}.
\]
\end{exam}

\section{Extension to Variable Intervals}
The methodology developed in this paper can be extended to handle the case where the interval $I$ itself depends on the random variable $x$, i.e.,
\[
G(t, x) = 0, \quad t \in I(x) = [l(x), u(x)],
\]
where $l(x), u(x) : \R^n \rightarrow \R$ are given functions.

This situation can be reduced to the fixed-interval case by restricting the density $\rho(x)$ to the set
\[
\Omega = \bigcup_{\hat{t}\in \R} \left( L(G(t,x), \{\hat{t}\}) \cap \{ x \in \R^n \mid l(x) \leq \hat{t} \leq u(x) \} \right).
\]
Then the expected number of roots in the variable interval $I(x)$ is given by
\[
\mathbb{E}N(G, \rho, I(\cdot)) = \mathbb{E}N(G, \rho \mathbb{I}(x, \Omega), \R).
\]

Similarly, we can define the intensity function for the variable interval case as
\[
T(\hat{t}) = \int_{L(G(t, x), \{\hat{t}\}) \cap \{ x \mid l(x) \leq \hat{t} \leq u(x) \}} \rho(x) \frac{|G'(\hat{t}, x)|}{\|\nabla G(\hat{t}, x)\|}  d\mathbb{H}^{n-1}.
\]
Then, under some suitable conditions analogous to those in the fixed-interval case, we could have
\[
\mathbb{E}N(G, \rho, I(\cdot)) = \int_{\R} T(t)  dt.
\]

This extension significantly broadens the applicability of our framework to problems where the domain of interest depends on the random parameters, such as in certain problems involving random domains or stopping times.

\section{Conclusion}
We have introduced a novel framework for computing the expected number of real roots of a stochastic function on a given interval. The classical Kac--Rice formula relies on the joint density of the function and its derivative, a requirement that severely limits its applicability. By studying the cumulative expectation function $\varphi_{G,\rho}(t)=\mathbb{E}N(G,\rho,[l,t])$ and establishing its absolute continuity, we derive two complementary computational schemes:
{Scheme 1} expresses the expectation as the derivative of a variable‑domain integral, under assumptions that are considerably weaker than those of the Kac--Rice formula;
{Scheme 2} provides an explicit integral representation free of joint densities and variable‑domain differentiation, under slightly stronger but easily verifiable conditions.

The framework is illustrated in detail for the linear case $G(t,x)=f(t)^\top x-f_0(t)$, yielding explicit formulas for Gaussian and uniform distributions and several new analytical results. The method also extends naturally to problems with variable intervals, widening its potential applications.

Our approach offers a rigorous and computationally more accessible alternative to traditional joint‑density‑based methods. Possible future directions include further relaxation of technical assumptions, extension to higher‑dimensional root‑counting problems, and the development of efficient numerical algorithms based on the obtained formulas.

In summary, this work provides a powerful new tool for analyzing random equations, with implications for stochastic analysis, statistical inference, and related fields.

\section*{Acknowledgments}
The author thanks the anonymous reviewers for their constructive comments and suggestions that helped improve the paper.

\section*{Funding}
This work was partially supported by the Jiangsu Provincial Scientific Research Center of Applied Mathematics [Grant No. BK20233002].

\begin{appendix}

\section{Proofs of Technical Results}
This appendix contains the detailed proofs of the lemmas and corollaries stated in the main text.

\subsection{Proof of Corollary~\ref{cor-1}}
\proof
For any $x \in L(\{P(s,x),Q(s,x)\},\hat{I})$, since $s = u(x)$ is the unique implicit function of $P(s,x) = 0$ on $\hat{I} \times \Omega$, we have $\eta(P(s,x),\hat{I},x)=\{u(x)\}$, and $Q(u(x),x)=0$. This implies $\{u(x)\} \subseteq \eta(Q(s,x),\hat{I},x)$, so $x \in \Xi(P(s,x), Q(s,x),\hat{I})$. Thus,
\[
L(\{P(s,x),Q(s,x)\},\hat{I}) \subseteq \Xi(P(s,x),Q(s,x),\hat{I}).
\]
Combining this with Property \ref{pro-1}(8), we obtain
\[
L(\{P(s,x),Q(s,x)\},\hat{I}) = \Xi(P(s,x),Q(s,x),\hat{I}).
\]

Since $s = u(x)$ is the unique implicit function of $P(s,x) = 0$ on $\hat{I} \times \Omega$, consider the set:
\[
\Delta(\varepsilon) = \left\{ x \in \Omega \middle| s = u(x),\ \varepsilon \leq \frac{|P'(s,x)|}{\|\nabla P(s,x)\|} \leq \frac{1}{\varepsilon} \right\}.
\]
Consider the integral
\[
\int_{\mathbb{R}^n} \mathbb{I}(x, \Delta(\varepsilon)) \rho(x) \frac{|P'(s,x)|}{\|\nabla P(s,x)\|} dx,
\]
where $s = u(x)$. Note that
\[
\mathbb{I}(x, \Delta(\varepsilon)) \rho(x) \frac{|P'(s,x)|}{\|\nabla P(s,x)\|} \leq \frac{\rho(x)}{\varepsilon}
\]
is integrable on $\mathbb{R}^n$. Moreover, on $\Delta(\varepsilon)$, since $\displaystyle\frac{\|\nabla P(s,x)\|}{|P'(s,x)|} \leq \frac{1}{\varepsilon}$, the function $s = u(x)$ satisfies the Lipschitz condition $\|\nabla u(x)\| = \displaystyle\frac{\|\nabla P(s,x)\|}{|P'(s,x)|} \leq \frac{1}{\varepsilon}$ on $\Delta(\varepsilon)$. We can extend $s = u(x)$ to the entire $\mathbb{R}^n$ while preserving the Lipschitz condition.

Using Lemma~\ref{coarea-1}, we get:
\[
\int_{\Delta(\varepsilon)} \rho(x) dx = \int_{\mathbb{R}^n} \mathbb{I}(x, \Delta(\varepsilon)) \rho(x) \frac{|P'(s,x)|}{\|\nabla P(s,x)\|} \|\nabla u(x)\| dx
\]
\[
= \int_{\mathbb{R}} d\hat{s} \int_{u^{-1}(\hat{s})} \mathbb{I}(x, \Delta(\varepsilon)) \rho(x) \frac{|P'(\hat{s},x)|}{\|\nabla P(\hat{s},x)\|} d\mathbb{H}^{n-1}
\]
\[
= \int_{u(\Delta(\varepsilon))} d\hat{s} \int_{L(P(s,x),\{\hat{s}\}) \cap \Delta(\varepsilon)} \rho(x) \frac{|P'(\hat{s},x)|}{\|\nabla P(\hat{s},x)\|} d\mathbb{H}^{n-1}.
\]
Note that
\[
\int_{\Delta(\varepsilon)} \rho(x) dx = \int_{\mathbb{R}^n} \mathbb{I}(x, \Delta(\varepsilon)) \rho(x) dx \leq \int_{\mathbb{R}^n} \rho(x) dx.
\]

Because $\Delta(\varepsilon)=\Omega \cap (L(P(s,x),\hat{I}) \setminus \Psi(\varepsilon))$, where
\[
\Psi(\varepsilon)=L(\{P(s,x),Q(s,x)\},\hat{I}) \cup \Psi_1(\varepsilon) \cup \Psi_2(\varepsilon),
\]
\[
\Psi_1(\varepsilon)=\left\{ x \in \Omega \middle| s = u(x),\ \frac{|P'(s,x)|}{\|\nabla P(s,x)\|} \leq \varepsilon \right\},
\]
\[
\Psi_2(\varepsilon)=\left\{ x \in \Omega \middle| s = u(x),\ \frac{|P'(s,x)|}{\|\nabla P(s,x)\|} \geq \frac{1}{\varepsilon} \right\},
\]
and clearly
\[
\lim_{\varepsilon \to 0_+} \Psi_1(\varepsilon) = L(\{P(s,x), P'(s,x)\},\hat{I}) \setminus L(\{P(s,x), \nabla P(s,x)\},\hat{I}),
\]
\[
\lim_{\varepsilon \to 0_+} \Psi_2(\varepsilon) = L(\{P(s,x), \nabla P(s,x)\},\hat{I}) \setminus L(\{P(s,x), P'(s,x)\},\hat{I}),
\]
we have
\[
\lim_{\varepsilon \to 0_+} \Psi(\varepsilon) = L(\{P(s,x),Q(s,x)\},\hat{I}) = \Xi(P(s,x),Q(s,x),\hat{I}) = \tau(\hat{I}),
\]
\[
\lim_{\varepsilon \to 0_+} \Delta(\varepsilon) = \Omega \cap (L(P(s,x),\hat{I}) \setminus \tau(\hat{I})) = \Omega \cap \tau(\hat{I}).
\]

We have
\[
\lim_{\varepsilon \to 0_+} \int_{\Delta(\varepsilon)} \rho(x) dx = \lim_{\varepsilon \to 0_+} \int_{\mathbb{R}^n} \mathbb{I}(x, \Delta(\varepsilon)) \rho(x) dx = \int_{\mathbb{R}^n} \lim_{\varepsilon \to 0_+} \mathbb{I}(x, \Delta(\varepsilon)) \rho(x) dx
\]
\[
= \int_{\mathbb{R}^n} \mathbb{I}(x, \Omega \setminus \tau(\hat{I})) \rho(x) dx = \int_{\Omega \setminus \tau(\hat{I})} \rho(x) dx.
\]
Because $\displaystyle\int_{\Omega \cap \tau(\hat{I})} \rho(x) dx = 0$, we have $\displaystyle\int_{\Omega \setminus \tau(\hat{I})} \rho(x) dx = \int_{\Omega} \rho(x) dx$.

Similarly, taking the limit on the right-hand side, and since $$T(\hat{s}) = \int_{L(P(s,x),\{\hat{s}\})} \rho(x) \frac{|P'(\hat{s},x)|}{\|\nabla P(\hat{s},x)\|} d\mathbb{H}^{n-1}$$ is integrable on $I$, we obtain
\[
\lim_{\varepsilon \to 0_+} \int_{u(\Delta(\varepsilon))} d\hat{s} \int_{L(P(s,x),\{\hat{s}\}) \cap \Delta(\varepsilon)} \rho(x) \frac{|P'(\hat{s},x)|}{\|\nabla P(\hat{s},x)\|} d\mathbb{H}^{n-1}
\]
\[
= \int_{\mathbb{R}} d\hat{s} \lim_{\varepsilon \to 0_+} \mathbb{I}(\hat{s}, u(\Delta(\varepsilon))) \int_{L(P(s,x),\{\hat{s}\}) \cap \Delta(\varepsilon)} \rho(x) \frac{|P'(\hat{s},x)|}{\|\nabla P(\hat{s},x)\|} d\mathbb{H}^{n-1}.
\]
Note that $\Delta(\varepsilon_1) \subseteq \Delta(\varepsilon_2)$ for $\varepsilon_1 \geq \varepsilon_2$, so $\displaystyle\lim_{\varepsilon \to 0_+} \Delta(\varepsilon) = \bigcup_{\varepsilon > 0} \Delta(\varepsilon) = \Omega \setminus \tau(\hat{I})$.

For any monotonically decreasing sequence $\varepsilon_k \to 0_+$, $\Delta(\varepsilon_k)$ is a monotonically increasing sequence of sets, so
\[
\lim_{k \to +\infty} u(\Delta(\varepsilon_k)) = \bigcup_{k > 0} u(\Delta(\varepsilon_k)).
\]
For any $y \in \bigcup_{k > 0} u(\Delta(\varepsilon_k))$, there exists $k^* > 0$ such that $y \in u(\Delta(\varepsilon_{k^*}))$, meaning there exists $x \in \Delta(\varepsilon_{k^*})$ with $u(x) = y$. Note that $x \in \bigcup_{k > 0} \Delta(\varepsilon_k)$, so $y \in u(\bigcup_{k > 0} \Delta(\varepsilon_k))$.

Conversely, for any $y \in u(\bigcup_{k > 0} \Delta(\varepsilon_k))$, there exists $x \in \bigcup_{k > 0} \Delta(\varepsilon_k)$, i.e., there exists $k^* > 0$ such that $x \in \Delta(\varepsilon_{k^*})$, so $y \in u(\Delta(\varepsilon_{k^*})) \subseteq \bigcup_{k > 0} u(\Delta(\varepsilon_k))$.

Thus, we have $$\lim_{\varepsilon \to 0_+} u(\Delta(\varepsilon)) = u(\bigcup_{\varepsilon > 0} \Delta(\varepsilon)) = u(\Omega \setminus \tau(\hat{I})).$$

Therefore,
\[
\int_{\mathbb{R}} d\hat{s} \lim_{\varepsilon \to 0_+} \mathbb{I}(\hat{s}, u(\Delta(\varepsilon))) \int_{L(P(s,x),\{\hat{s}\}) \cap \Delta(\varepsilon)} \rho(x) \frac{|P'(\hat{s},x)|}{\|\nabla P(\hat{s},x)\|} d\mathbb{H}^{n-1}
\]
\[
= \int_{\mathbb{R}} d\hat{s} \mathbb{I}(\hat{s}, u(\Omega \setminus \tau(\hat{I}))) \lim_{\varepsilon \to 0_+} \int_{L(P(s,x),\{\hat{s}\}) \cap \Delta(\varepsilon)} \rho(x) \frac{|P'(\hat{s},x)|}{\|\nabla P(\hat{s},x)\|} d\mathbb{H}^{n-1}
\]
\[
= \int_{\mathbb{R}} d\hat{s} \mathbb{I}(\hat{s}, u(\Omega \setminus \tau(\hat{I}))) \int_{\Omega \cap L(P(s,x),\{\hat{s}\}) \setminus \tau(\{\hat{s}\})} \rho(x) \frac{|P'(\hat{s},x)|}{\|\nabla P(\hat{s},x)\|} d\mathbb{H}^{n-1}
\]

Since $T(\hat{s})$ is Lebesgue integrable on $\hat{I}$, it is finite almost everywhere on $\hat{I}$.
\[
= \int_{u(\Omega \setminus \tau(\hat{I}))} d\hat{s} \int_{\Omega \cap L(P(s,x),\{\hat{s}\}) \setminus \tau(\{\hat{s}\})} \rho(x) \frac{|P'(\hat{s},x)|}{\|\nabla P(\hat{s},x)\|} d\mathbb{H}^{n-1},
\]
\[
= \int_{u(\Omega \setminus \tau(\hat{I}))} d\hat{s} \int_{\Omega \cap L(P(s,x),\{\hat{s}\})} \rho(x) \frac{|P'(\hat{s},x)|}{\|\nabla P(\hat{s},x)\|} d\mathbb{H}^{n-1} \]\[- \int_{u(\Omega \setminus \tau(\hat{I}))} d\hat{s} \int_{\Omega \cap L(P(s,x),\{\hat{s}\}) \cap \tau(\{\hat{s}\})} \rho(x) \frac{|P'(\hat{s},x)|}{\|\nabla P(\hat{s},x)\|} d\mathbb{H}^{n-1}.
\]

Because for almost every $s \in \hat{I}$,
\[
\int_{\Omega \cap L(P(s,x),\{\hat{s}\}) \cap \tau(\{\hat{s}\})} \rho(x) \frac{|P'(\hat{s},x)|}{\|\nabla P(\hat{s},x)\|} d\mathbb{H}^{n-1} = \int_{\Omega \cap \tau(\{\hat{s}\})} \rho(x) \frac{|P'(\hat{s},x)|}{\|\nabla P(\hat{s},x)\|} d\mathbb{H}^{n-1} = 0,
\]
so
\[
\int_{u(\Omega \setminus \tau(\hat{I}))} d\hat{s} \int_{\Omega \cap L(P(s,x),\{\hat{s}\}) \setminus \tau(\{\hat{s}\})} \rho(x) \frac{|P'(\hat{s},x)|}{\|\nabla P(\hat{s},x)\|} d\mathbb{H}^{n-1}
\]
\[
= \int_{u(\Omega \setminus \tau(\hat{I}))} d\hat{s} \int_{\Omega \cap L(P(s,x),\{\hat{s}\})} \rho(x) \frac{|P'(\hat{s},x)|}{\|\nabla P(\hat{s},x)\|} d\mathbb{H}^{n-1}.
\]

Moreover, because $u(\Omega \setminus \tau(\hat{I})) = \hat{I}$, we finally obtain
\[
\int_{\Omega} \rho(x) dx = \int_{\hat{I}} d\hat{s} \int_{\Omega \cap L(P(s,x),\{\hat{s}\})} \rho(x) \frac{|P'(\hat{s},x)|}{\|\nabla P(\hat{s},x)\|} d\mathbb{H}^{n-1}.
\]
\qed

\subsection{Proof of Corollary~\ref{cor-2}}
\proof
Given $M > 0$, $\varepsilon > 0$, consider the set:
\[
\Delta(M,\varepsilon) = \left\{ x \in B(0_n, M) \cap L(P(s,x),\hat{I}) \,\middle|\, \exists s \in \hat{I},\ P(s,x) = 0,\ \varepsilon \leq \frac{|P'(s,x)|}{\|\nabla P(s,x)\|} \leq \frac{1}{\varepsilon} \right\}.
\]
For any $x^* \in \Delta(M,\varepsilon)$, there exists $s^* \in \hat{I}$ such that: there exists a neighborhood of $(s^*, x^*)$
\[
(B(s^*,\delta) \times B(x^*,\delta)) \cap (\hat{I} \times B(0_n,M)),
\]
on which there exists a unique implicit function denoted $s = u^*(x) : \R^n \rightarrow \R$.

For any $x_k \in L(P(s,x),\hat{I}), x_k \rightarrow x^*$, $\exists s_k \in \hat{I}$ satisfying $P(s_k,x_k) = 0$. Since $\hat{I}$ is compact, $\exists s^*, s_{k_i} \rightarrow s^*$, and by the continuity of $P(s,x)$, we have $P(s^*,x^*) = 0$, i.e., $x^* \in L(P(s,x),\hat{I})$. This shows that $L(P(s,x),\hat{I})$ is a closed set, hence $\Delta(M,\varepsilon)$ is compact.

By the Heine–Borel theorem, we obtain a finite number of regions
\[
\Theta_i(M,\varepsilon) \subseteq B(s_i,\delta_i) \cap \hat{I}, \quad \Sigma_i(M,\varepsilon) \subseteq B(x_i,\delta_i) \cap \Delta(M,\varepsilon)
\]
\[
\Sigma_i(M,\varepsilon) \cap \Sigma_j(M,\varepsilon) = \emptyset \ (i \neq j), \quad \bigcup_i \Sigma_i(M,\varepsilon) = \Delta(M,\varepsilon).
\]

Then
\[
\mathbb{I}(x, \Sigma_i(M,\varepsilon)) \rho(x) \frac{|P'(s,x)|}{\|\nabla P(s,x)\|} \leq \frac{\rho(x)}{\varepsilon},
\]
\[
\|\nabla u_i(x)\| = \frac{\|\nabla P(s,x)\|}{|P'(s,x)|} \leq \frac{1}{\varepsilon},
\]
thus on $\Sigma_i(M,\varepsilon)$ the conditions for the Coarea Formula are satisfied, yielding
\[
\int_{\Sigma_i(M,\varepsilon)} \rho(x)  dx = \int_{\R^n} \mathbb{I}(x, \Sigma_i(M,\varepsilon)) \rho(x) \frac{|P'(s,x)|}{\|\nabla P(s,x)\|} \|\nabla u_i(x)\|  dx
\]
\[
= \int_{u_i(\Sigma_i(M,\varepsilon))} d\hat{s} \int_{L(P(s,x),\{\hat{s}\}) \cap \Sigma_i(M,\varepsilon)} \rho(x) \frac{|P'(\hat{s},x)|}{\|\nabla P(\hat{s},x)\|}  d\mathbb{H}^{n-1}
\]
\[
\leq \int_{\hat{I}} d\hat{s} \int_{L(P(s,x),\{\hat{s}\}) \cap \Sigma_i(M,\varepsilon)} \rho(x) \frac{|P'(\hat{s},x)|}{\|\nabla P(\hat{s},x)\|}  d\mathbb{H}^{n-1}.
\]

Therefore,
\[
\sum_{i} \int_{\Sigma_i(M,\varepsilon)} \rho(x)  dx \leq \sum_i \int_{\hat{I}} d\hat{s} \int_{L(P(s,x),\{\hat{s}\}) \cap \Sigma_i(M,\varepsilon)} \rho(x) \frac{|P'(\hat{s},x)|}{\|\nabla P(\hat{s},x)\|}  d\mathbb{H}^{n-1},
\]
i.e.,
\[
\int_{\bigcup_i \Sigma_i(M,\varepsilon)} \rho(x)  dx \leq \int_{\hat{I}} d\hat{s} \int_{L(P(s,x),\{\hat{s}\}) \cap (\bigcup_i \Sigma_i(M,\varepsilon))} \rho(x) \frac{|P'(\hat{s},x)|}{\|\nabla P(\hat{s},x)\|}  d\mathbb{H}^{n-1}
\]
\[
\leq \int_{\hat{I}} d\hat{s} \int_{L(P(s,x),\{\hat{s}\})} \rho(x) \frac{|P'(\hat{s},x)|}{\|\nabla P(\hat{s},x)\|}  d\mathbb{H}^{n-1}.
\]

Note that
\[
\lim_{\varepsilon \to 0_+} \bigcup_i \Sigma_i(M,\varepsilon) = \lim_{\varepsilon \to 0_+} \Delta(M,\varepsilon) = B(0_n, M) \cap (L(P(s,x),\hat{I}) \setminus \Pi),
\]
where
\[
\Pi = \{ x \in L(P(s,x),\hat{I}) \mid \forall s \in \{s\in \hat{I} | \ P(s,x) = 0\}, Q(s,x) = 0 \}=\Xi(P(s,x), Q(s,x), \hat{I}).
\]

By the condition $\displaystyle\int_{\Xi(P(s,x), Q(s,x), \hat{I})} \rho(x)  dx = 0$,
we finally obtain
\[
\lim_{M \to +\infty} \lim_{\varepsilon \to 0_+} \int_{\bigcup_i \Sigma_i(M,\varepsilon)} \rho(x)  dx = \int_{L(P(s,x),\hat{I}) \setminus \Pi} \rho(x)  dx = \int_{L(P(s,x),\hat{I})} \rho(x)  dx
\]
\[
\leq \int_{\hat{I}} d\hat{s} \int_{L(P(s,x),\{\hat{s}\})} \rho(x) \frac{|P'(\hat{s},x)|}{\|\nabla P(\hat{s},x)\|}  d\mathbb{H}^{n-1}.
\]
\qed

\subsection{Proof of Corollary~\ref{cor-3}}
\proof This proof is similar as Corollary~\ref{cor-2}.
Given $M > 0$, $\varepsilon > 0$, consider the set:
\[
\Delta(M,\varepsilon) = \left\{ x \in B(0_n, M) \cap L(P(s,x),\hat{I}^2) \,\middle|\, \exists s \in \hat{I}^2,\ P(s,x) = 0_2,\ \right. 
\]
\[
\frac{|P'_1(s_1,x)| |P'_2(s_2,x)|}{\sqrt{V(s,x)}} \leq \frac{1}{\varepsilon},
\left.\|V_2(s,x)^{-1} V_1(s,x) V_2(s,x)^{-1}\| \leq \frac{1}{\varepsilon} \right\},
\]
where
\[
V_2(s,x) = \begin{pmatrix} 
P_1'(s_1,x) & 0 \\
0 & P_2'(s_2,x)
\end{pmatrix}, \]\[
V_1(s,x) = \begin{pmatrix}  
\nabla P_{1}(s_1,x)^\top \nabla P_1(s_1,x) & \nabla P_{1}(s_1,x)^\top \nabla P_2(s_2,x) \\
\nabla P_{2}(s_2,x)^\top \nabla P_1(s_1,x) & \nabla P_{2}(s_2,x)^\top \nabla P_2(s_2,x)
\end{pmatrix},
\]
\[
V(s,x) = \det(V_1(s,x)).
\]

For any $x^* \in \Delta(M,\varepsilon)$, there exists $s^* \in \hat{I}^2$ such that: there exists a neighborhood of $(s^*, x^*)$
\[
(B(s^*,\delta) \times B(x^*,\delta)) \cap (\hat{I}^2 \times B(0_n,M)),
\]
on which there exists a unique implicit function denoted $s = u^*(x) : \R^n \rightarrow \R^2$.

For any $x_k \in L(P(s,x),\hat{I}^2), x_k \rightarrow x^*$, $\exists s_k \in \hat{I}^2$ satisfying $P(s_k,x_k) = 0_2$. Since $\hat{I}^2$ is compact, $\exists s^*, s_{k_i} \rightarrow s^*$, and by the continuity of $P(s,x)$, we have $P(s^*,x^*) = 0_2$, i.e., $x^* \in L(P(s,x),\hat{I}^2)$. This shows that $L(P(s,x),\hat{I}^2)$ is a closed set, hence $\Delta(M,\varepsilon)$ is compact.

By the Heine–Borel theorem, we obtain a finite number of regions
\[
\Theta_i(M,\varepsilon) \subseteq B(s_i,\delta_i) \cap \hat{I}^2, \quad \Sigma_i(M,\varepsilon) \subseteq B(x_i,\delta_i) \cap \Delta(M,\varepsilon)
\]
\[
\Sigma_i(M,\varepsilon) \cap \Sigma_j(M,\varepsilon) = \emptyset \ (i \neq j), \quad \bigcup_i \Sigma_i(M,\varepsilon) = \Delta(M,\varepsilon).
\]

Then
\[
\left| \mathbb{I}(x, \Sigma_i(M,\varepsilon)) \rho(x) \frac{|P_1'(s_1,x)| |P_2'(s_2,x)|}{\sqrt{V(s,x)}} \right| \leq \frac{\rho(x)}{\varepsilon},
\]
\[
\|\nabla u_i(x)\|^2 = \|V_2(s,x)^{-1} V_1(s,x) V_2(s,x)^{-1}\| \leq \frac{1}{\varepsilon}.
\]

Thus on $\Sigma_i(M,\varepsilon)$ the conditions for the Coarea Formula are satisfied, yielding
\[
\int_{\Sigma_i(M,\varepsilon)} \rho(x)  dx = \int_{\R^n} \mathbb{I}(x, \Sigma_i(M,\varepsilon)) \rho(x) \frac{|P_1'(s_1,x)| |P_2'(s_2,x)|}{\sqrt{V(s,x)}} \sqrt{\det(\nabla u_i(x)^\top \nabla u_i(x))}  dx 
\]
\[
= \int_{u_i(\Sigma_i(M,\varepsilon))} d\hat{s} \int_{u^{-1}_i(\hat{s}) \cap \Sigma_i(M,\varepsilon)} \rho(x) \frac{|P_1'(\hat{s}_1,x)| |P_2'(\hat{s}_2,x)|}{\sqrt{V(\hat{s},x)}}  d\mathbb{H}^{n-2}
\]
\[
= \int_{u_i(\Sigma_i(M,\varepsilon))} d\hat{s} \int_{L(P(s,x),\{\hat{s}\}) \cap \Sigma_i(M,\varepsilon)} \rho(x) \frac{|P_1'(\hat{s}_1,x)| |P_2'(\hat{s}_2,x)|}{\sqrt{V(\hat{s},x)}}  d\mathbb{H}^{n-2}
\]
\[
\leq \int_{\hat{I}^2} d\hat{s} \int_{L(P(s,x),\{\hat{s}\}) \cap \Sigma_i(M,\varepsilon)} \rho(x) \frac{|P_1'(\hat{s}_1,x)| |P_2'(\hat{s}_2,x)|}{\sqrt{V(\hat{s},x)}}  d\mathbb{H}^{n-2}.
\]

Therefore,
\[
\sum_{i} \int_{\Sigma_i(M,\varepsilon)} \rho(x)  dx = \int_{\bigcup_i \Sigma_i(M,\varepsilon)} \rho(x)  dx
\]
\[
\leq \sum_i \int_{\hat{I}^2} d\hat{s} \int_{L(P(s,x),\{\hat{s}\}) \cap \Sigma_i(M,\varepsilon)} \rho(x) \frac{|P_1'(\hat{s}_1,x)| |P_2'(\hat{s}_2,x)|}{\sqrt{V(\hat{s},x)}}  d\mathbb{H}^{n-2}
\]
\[
= \int_{\hat{I}^2} d\hat{s} \int_{L(P(s,x),\{\hat{s}\}) \cap (\bigcup_i \Sigma_i(M,\varepsilon))} \rho(x) \frac{|P_1'(\hat{s}_1,x)| |P_2'(\hat{s}_2,x)|}{\sqrt{V(\hat{s},x)}}  d\mathbb{H}^{n-2}
\]
\[
\leq \int_{\hat{I}^2} d\hat{s} \int_{L(P(s,x),\{\hat{s}\})} \rho(x) \frac{|P_1'(\hat{s}_1,x)| |P_2'(\hat{s}_2,x)|}{\sqrt{V(\hat{s},x)}}  d\mathbb{H}^{n-2}.
\]

Now note that
\[
\lim_{\varepsilon \to 0_+} \bigcup_i \Sigma_i(M,\varepsilon) = \lim_{\varepsilon \to 0_+} \Delta(M,\varepsilon) = B(0_n, M) \cap (L(P(s,x),\hat{I}^2) \setminus \Pi)
\]
where
\[
\Pi = \{ x \in L(P(s,x),\hat{I}^2) \mid \forall s \in \{s \in \hat{I}^2 | P(s,x) = 0_2 \}, Q(s,x) = 0 \} \]\[= \Xi(P(s,x), Q(s,x),\hat{I}^2).
\]

By the condition 
\[
\int_{\Xi(P(s,x), Q(s,x),\hat{I}^2)} \rho(x)  dx = 0,
\]
we finally obtain
\[
\lim_{M \to +\infty} \lim_{\varepsilon \to 0_+} \int_{\bigcup_i \Sigma_i(M,\varepsilon)} \rho(x)  dx = \int_{L(P(s,x),\hat{I}^2) \setminus \Pi} \rho(x)  dx = \int_{L(P(s,x),\hat{I}^2)} \rho(x)  dx
\]
\[
\leq \int_{\hat{I}^2} d\hat{s} \int_{L(P(s,x),\{\hat{s}\})} \rho(x) \frac{|P_1'(\hat{s}_1,x)| |P_2'(\hat{s}_2,x)|}{\sqrt{V(\hat{s},x)}}  d\mathbb{H}^{n-2}.
\]
\qed

\subsection{Proof of Corollary~\ref{cor-4}}
\proof
Given $M > 0$, $\varepsilon > 0$, consider the set:
\[
\Delta(M,\varepsilon) = \left\{ x \in B(0_n, M) \cap L(P_2(s,x),\hat{I}) \,\middle|\, \exists s \in \hat{I},\ P_2(s,x) = 0,\ 
\|\nabla P_1(x)\| \geq \varepsilon,\ \right.
\]
\[
\left. \frac{\|\nabla P_2(s,x)\|}{|P_2'(s,x)|} \leq \frac{1}{\varepsilon},
\frac{|P_2'(s,x)| \|\nabla P_1(x)\|}{\sqrt{V(s,x)}} \leq \frac{1}{\varepsilon} \right\}.
\]

For any $x^* \in \Delta(M,\varepsilon)$, there exists $s^* \in \hat{I}$ such that: there exists a neighborhood of $(s^*, x^*)$
\[
(B(s^*,\delta) \times B(x^*,\delta)) \cap (\hat{I} \times B(0_n, M)),
\]
on which there exists a unique implicit function denoted $s = u^*(x) : \R^n \rightarrow \R$.

For any $x_k \in L(P_2(s,x),\hat{I}), x_k \rightarrow x^*$, $\exists s_k \in \hat{I}$ satisfying $P_2(s_k,x_k) = 0$. Since $\hat{I}$ is compact, $\exists s^*, s_{k_i} \rightarrow s^*$, and by the continuity of $P_2(s,x)$, we have $P_2(s^*,x^*) = 0$, i.e., $x^* \in L(P_2(s,x),\hat{I})$. This shows that $L(P_2(s,x),\hat{I})$ is a closed set, hence $\Delta(M,\varepsilon)$ is compact.

By the Heine–Borel theorem, we obtain a finite number of regions
\[
\Theta_i(M,\varepsilon) \subseteq B(s_i,\delta_i) \cap \hat{I}, \quad \Sigma_i(M,\varepsilon) \subseteq B(x_i,\delta_i) \cap \Delta(M,\varepsilon)
\]
\[
\Sigma_i(M,\varepsilon) \cap \Sigma_j(M,\varepsilon) = \emptyset \ (i \neq j), \quad \bigcup_i \Sigma_i(M,\varepsilon) = \Delta(M,\varepsilon).
\]

Then
\[
\|\nabla u_i(x)\| = \frac{\|\nabla P_2(s,x)\|}{|P_2'(s,x)|} \leq \frac{1}{\varepsilon},
\]
\[
\left| \mathbb{I}(x, \Omega \cap \Sigma_i(M,\varepsilon)) \rho(x) \frac{|P_2'(\hat{s},x)| \|\nabla P_1(x)\|}{\sqrt{V(\hat{s},x)}} \right| \leq \frac{\rho(x)}{\varepsilon}.
\]

According to Lemma~\ref{coarea-2} (the Coarea formula on rectifiable sets), we have:
\[
\int_{\Omega \cap \Sigma_i(M,\varepsilon)} \rho(x)  d\mathbb{H}^{n-1}\]\[ = \int_{\Omega \cap B(0_n, M)} \mathbb{I}(x, \Omega \cap \Sigma_i(M,\varepsilon)) \rho(x) \frac{|P_2'(s,x)| \|\nabla P_1(x)\|}{\sqrt{V(s,x)}} \frac{\sqrt{V(s,x)}}{|P_2'(s,x)| \|\nabla P_1(x)\|}  d\mathbb{H}^{n-1}
\]
\[
= \int_{u_i(\Omega \cap \Sigma_i(M,\varepsilon))} d\hat{s} \int_{\Omega \cap (u_i^{-1}(\hat{s}) \cap \Sigma_i(M,\varepsilon))} \rho(x) \frac{|P_2'(\hat{s},x)| \|\nabla P_1(x)\|}{\sqrt{V(\hat{s},x)}}  d\mathbb{H}^{n-2}
\]
\[
= \int_{u_i(\Omega \cap \Sigma_i(M,\varepsilon))} d\hat{s} \int_{\Omega \cap (L(P_2(s,x),\{\hat{s}\}) \cap \Sigma_i(M,\varepsilon))} \rho(x) \frac{|P_2'(\hat{s},x)| \|\nabla P_1(x)\|}{\sqrt{V(\hat{s},x)}}  d\mathbb{H}^{n-2}
\]
\[
\leq \int_{\hat{I}} d\hat{s} \int_{\Omega \cap (L(P_2(s,x),\{\hat{s}\}) \cap \Sigma_i(M,\varepsilon))} \rho(x) \frac{|P_2'(\hat{s},x)| \|\nabla P_1(x)\|}{\sqrt{V(\hat{s},x)}}  d\mathbb{H}^{n-2}.
\]

Therefore,
\[
\sum_{i} \int_{\Omega \cap \Sigma_i(M,\varepsilon)} \rho(x)  d\mathbb{H}^{n-1} = \int_{\Omega \cap \bigcup_i \Sigma_i(M,\varepsilon)} \rho(x)  d\mathbb{H}^{n-1}
\]
\[
\leq \sum_i \int_{\hat{I}} d\hat{s} \int_{\Omega \cap (L(P_2(s,x),\{\hat{s}\}) \cap \Sigma_i(M,\varepsilon))} \rho(x) \frac{|P_2'(\hat{s},x)| \|\nabla P_1(x)\|}{\sqrt{V(\hat{s},x)}}  d\mathbb{H}^{n-2}
\]
\[
= \int_{\hat{I}} d\hat{s} \int_{\Omega \cap (L(P_2(s,x),\{\hat{s}\}) \cap (\bigcup_i \Sigma_i(M,\varepsilon)))} \rho(x) \frac{|P_2'(\hat{s},x)| \|\nabla P_1(x)\|}{\sqrt{V(\hat{s},x)}}  d\mathbb{H}^{n-2}
\]
\[
\leq \int_{\hat{I}} d\hat{s} \int_{\Omega \cap L(P_2(s,x),\{\hat{s}\})} \rho(x) \frac{|P_2'(\hat{s},x)| \|\nabla P_1(x)\|}{\sqrt{V(\hat{s},x)}}  d\mathbb{H}^{n-2}.
\]

Now note that
\[
\lim_{\varepsilon \to 0_+} \Omega \cap \bigcup_i \Sigma_i(M,\varepsilon) = \lim_{\varepsilon \to 0_+} \Omega \cap \Delta(M,\varepsilon) = B(0_n, M) \cap \Omega \cap (L(P_2(s,x),\hat{I}) \setminus \Pi)
\]
where
\[
\Pi = \{ x \in L(P_2(s,x),\hat{I}) \mid \forall s \in \{s \in \hat{I} | \ P_2(s,x) = 0 \}, P'_2(s,x) \|\nabla P_1(x)\| V(s,x) = 0 \}\]\[ = \Xi(P_2(s, x), Q(s, x), \hat{I}).
\]

By the condition 
\[
\int_{\Omega \cap \Xi(P_2(s, x), Q(s, x), \hat{I})} \rho(x)  d\mathbb{H}^{n-1} = 0,
\]
\[
\lim_{M \to +\infty} \lim_{\varepsilon \to 0_+} \int_{\Omega \cap \bigcup_i \Sigma_i(M,\varepsilon)} \rho(x)  d\mathbb{H}^{n-1} = \int_{\Omega \cap (L(P_2(s,x),\hat{I}) \setminus \Pi)} \rho(x)  d\mathbb{H}^{n-1}
\]
\[
= \int_{\Omega \cap L(P_2(s,x),\hat{I})} \rho(x)  d\mathbb{H}^{n-1} \leq \int_{\hat{I}} d\hat{s} \int_{\Omega \cap L(P_2(s,x),\{\hat{s}\})} \rho(x) \frac{|P_2'(\hat{s},x)| \|\nabla P_1(x)\|}{\sqrt{V(\hat{s},x)}}  d\mathbb{H}^{n-2}.
\]
\qed

\end{appendix}
%\linenumbers %% 如果需要行号可以取消注释
\end{document}